\documentclass{conm-p-l}
\usepackage{amsmath,m-pictex}
\usepackage{enumerate}
\usepackage{amssymb,geometry,latexsym,color}
\usepackage{citesort}
\usepackage{xypic}
\xyoption{2cell}\UseTwocells
\UseHalfTwocells
\author{ Timothy Porter}
\address{School of Informatics, University of Wales Bangor, Bangor, Gwynedd, LL57 1UT, United Knigdom.}
\email{t.porter@bangor.ac.uk}

\title{Formal Homotopy Quantum Field Theories, II : Simplicial Formal Maps}
\newtheorem{theorem}{Theorem} 
\newtheorem{proposition}{Proposition} 
 
\newtheorem{lemma}{Lemma} 
  
\begin{document}

\begin{abstract}
Simplicial formal maps were introduced in the first paper of this series as a tool for studying Homotopy Quantum Field Theories with background a general homotopy 2-type. Here we continue their study, showing how a natural generalisation can handle much more general backgrounds.  The question of the geometric interpretation of these formal maps is partially answered in terms of combinatorial bundles.  This suggests new interpretations of HQFTs.

{\noindent\bf A. M. S. Classification:}  Primary:  18G50. Secondary: 55P99, 57R56, 81T45. 
\end{abstract}
\maketitle

\section*{Introduction}
 \medskip
In the Homotopy Quantum Field Theories introduced in \cite{turaev:hqft1,turaev:hqft2}, an important role is played by the background space, here denoted  $B$.   The objects of study are manifolds with extra structure and that extra structure is given by a `characteristic map' from the manifold to the target background space $B$.  These `$B$-manifolds' and $B$-cobordisms are then studied using tools similar to those of Topological Quantum Field Theories. In those initial papers, one axiom in the theory was unnecessarily strong and resulted in eliminating structure in $B$ above its $d$-type, when the manifolds concerned were of dimension $d$. A modified version with change to one axiom (see below) was introduced by Rodrigues, \cite{rodrigues}. This gave dependence of $(d+1)$-HQFTs over $B$ on the $(d+1)$-type of $B$. This was used by  Brightwell and Turner, \cite{B&T}, and Turner and Willerton, \cite{t&w02}, to look at  (1+1)-HQFTs with background space a simply connected space.  Thus the results of \cite{turaev:hqft1} had classified (1+1)-HQFTs with background spaces which were 1-types and the more recent results handled simply connected spaces, classification results there being in terms of the second homotopy group of $B$.  It was therefore natural to try to classify such HQFTs for which the background space is a 2-type, a situation that would include both the previous cases.

In trying out ideas for adapting the existing theory, it seemed that (i) part of our basic theory seemed to work just as well even if we did not restrict to (1+1)-HQFTs, and (ii) for $(d+1)$-HQFTs, we could assume $B$ was the classifying space of a crossed complex, in the sense \cite{B&H1991}. Some of the methods worked in even greater generality namely when $B$ was the classifying space of a $(d+1)$-truncated simplicial group, and thus was a general $(d+1)$-type.  This led us to a concept of \emph{simplicial formal map}, which provides an algebraic / combinatorial model for the characteristic map $g : M\to B$ that specifies the basic background structure for the manifold $M$.   

We introduced formal $\mathcal{C}$-maps on 1- and 2-dimensional manifolds in \cite{TPVT:hqft1} for $\mathcal{C}$ a crossed module. Adapting the axioms of HQFTs to work with formal $\mathcal{C}$-maps rather than $B$-manifolds gave the notion of formal HQFT for those dimensions and gave some classification results in that setting.

In this paper we introduce another approach to these formal maps in more generality, which suggests additional geometric interpretations of them even in the dimensions for which they were originally introduced.  This should lead on to a state sum type approach to constructing more general formal HQFTs.  It will also link in these formal HQFTs to various other areas of the interface between mathematics and physics.

The key to our approach is to have good algebraic models for homotopy types. We have restricted attention, for the detailed development here, to crossed complexes. These are excellently structured algebraic / categorical models for certain homotopy types of spaces.  Our methods at present work best with such crossed complexes, work moderately well with simplicial groups and we do not at all know, or at least, not yet, how to handle the more difficult, but more interesting, weak $n$-groupoid models that have appeared  in the literature.

So as not to end up with too long a paper, we will assume that the reader has at hand the introductory papers \cite{turaev:hqft1,turaev:hqft2} and the thesis by Rodrigues, \cite{rodrigues}. (A summary of this theory is given in an introductory section.)

\textsc{Acknowledgements.}
This work was partially supported by a grant, GR/S17635/01, from the EPSRC, for a visit by Turaev to Bangor and Gregynog Hall.   
This paper continues the earlier joint one with Turaev, and would not have been written without his collaboration on the overall project. I would also like to acknowledge the help given by Ronnie Brown, who participated in many of the discussions both at Bangor and at Gregynog. His wealth of ideas, perspective and knowledge on crossed modules, crossed complexes and all the general `crossed menagery' were invaluable. This paper is dedicated both to Ross whose 60th birthday is this year and to Ronnie who beats him by 10 years. Other significant anniversaries are Max Kelly's at 75, and the 100th anniversary of Ehresmann. Each of these have made important contributions to the parts of higher dimensional algebra on which this work is based.

\tableofcontents


\section{Homotopy Quantum Field Theories}
For the convenience of the reader we have included here a brief introduction to HQFTs in general.

Fix an integer $n \geq 0$ and a field, $\mathbb{K}$. All vector spaces will be tacitly assumed to be finite dimensional.  Usually $\mathbb{K}$ can be replaced by a commutative ring merely by replacing finite dimensional vector spaces by projective $\mathbb{K}$-modules of finite type, but we will not do this here. 

\subsection{The category of $B$-manifolds and $B$-cobordisms}
The basic objects on which a $(d+1)$-homotopy quantum field theory is built are compact, oriented $d$-manifolds together with  maps to a `background' space, $B$.  This space $B$ will  be path connected with a fixed base point, $\ast$.  More precisely:

\medskip

\textsc{Definition.} A \emph{$B$-manifold} is a pair $(X, g)$, where $X$ is a closed oriented $d$-manifold (with a choice of base point $m_i$ in each connected component $X_i$ of $X$), and $g$ is a continuous map $g : X \to B$, called the \emph{characteristic map}, such that $g(m_i) = \ast$ for each base point $m_i$.

A \emph{$B$-isomorphism} between $B$-manifolds, $\phi : ( X, g) \to ( Y, h)$ is an isomorphism $\phi : X \to Y$ of the manifolds, preserving the orientation, taking base points into base points and such that $h\phi = g$.

\medskip

\textsc{Remark.} 
It may sometimes be the case that the manifolds under consideration will be differentiable and then `isomorphism' is interpreted as `diffeomorphism', but equally well we can position the theory in the category of PL-manifolds or triangulable topological manifolds with the obvious changes.  In fact for some of the time it is convenient to  develop constructions for simplicial complexes rather than manifolds, as it is triangulations that provide the basis for the combinatorial descriptions of the structures that we will be using.  

\medskip

Denote by $\mathbf{Man}(d,B)$ the category of $d$-dimensional $B$-manifolds and $B$-isomorphisms.  We define a `sum' operation on this category using disjoint union.  The disjoint union of $B$-manifolds is defined by
$$( X, g) \amalg ( Y, h) := ( X\amalg Y, g\amalg h),$$
with the obvious characteristic map, $g\amalg h : X \amalg Y \to B$.  With this `sum' operation, $\mathbf{Man}(n,B)$ becomes a symmetric monoidal category with the unit being given by the empty $B$-manifold, $\emptyset$, with the empty characteristic map.  Of course, this is an $d$-manifold by default.

\medskip

It is important to remember that $(X, g) \amalg \emptyset$  is not really the same as $(X,g)$, but is naturally isomorphic to it via the obvious $B$-isomorphism $$l_{( X, g)} : ( X, g) \amalg \emptyset\to ( X, g).$$Of course there is a similar $B$-isomorphism, $r_{(X, g)} : \emptyset\amalg(X, g)  \to (X, g).$  Likewise $(X,g)\amalg(Y,h)$ is a categorical coproduct so is only determined up to natural isomorphism.  These are, of course, the  problems in most naturally arising monoidal structures such as the monoidal category $(Vect, \otimes)$ of finite dimensional vector spaces with tensor product as the monoidal structure and they motivate and guide the theory of such structures.

\medskip

\textsc{Definition.}
A \emph{cobordism} $W : X_0 \to X_1$ is a compact oriented $(d + 1)$-manifold, $W$, whose boundary is the disjoint union of pointed  closed oriented  $d$-manifolds, $X_0$ and $X_1$, such that the orientation of $X_1$ (resp. $X_0$) is induced by that on $W$ (resp., is opposite to the one induced from that on $W$).  (The manifold $W$ is not considered as being pointed.) It may be convenient to write $\partial W = -X_0 \amalg X_1$ and also $\partial_- W = X_0 $ and $\partial_+ W = X_1$. We may refer to $X_0$ as the \emph{incoming boundary} and $X_1$ as the \emph{outgoing boundary}, in the usual way.

A $B$\emph{-cobordism}, $(W,F)$, from $(X_0,g)$ to $(X_1,h)$ is a cobordism $W : X_0 \to X_1$ endowed with a homotopy class of maps $F : W \to X$ relative to the boundary such that $F|_{X_0} = 
g$ and $F|_{X_1} = h$. (Generally we will not make a notational distinction between the homotopy class $F$ and any of its representatives.)  Finally a $B$\emph{-isomorphism of $B$-cobordisms}, $\psi : (W,F) \to (W^\prime, F^\prime)$, is an isomorphism $\psi : W \to W^\prime$ such that 
$$\psi (\partial_+W) = \partial_+W^\prime,$$
$$\psi (\partial_-W) = \partial_-W^\prime,$$
and $F^\prime \psi = F$, in the obvious sense of homotopy classes relative to the boundary.

We can glue $B$-cobordisms along their boundaries, or more generally, along a $B$-isomorphism between their boundaries, in the usual way, see Turaev, \cite{turaev:hqft1}, or Rodrigues, \cite{rodrigues}. For each $B$-manifold, $(X, g)$, there is a $B$-cobordism $(I\times X, 1_g) : (X,g) \to (X,g)$ with $1_g(t,x) = g(x)$ and where, as usual, $I$ denotes the unit interval. This cobordism will be called the \emph{identity cobordism} on $(X,g)$ and will be denoted $1_{(X,g)}$.

As for disjoint union of $B$-manifolds, we can define a disjoint union of $B$-cobordisms, in the obvious way.
\medskip

\textsc{Remark.}
The detailed structure of $B$-cobordisms and the resulting category $\mathbf{HCobord}(d,B)$ is given in the Appendix to \cite{rodrigues} at least in the important case of differentiable $B$-manifolds.  This category is, technically, a monoidal category with strict duals and a homotopy quantum field theory will be  a symmetric monoidal functor from $\mathbf{HCobord}(d,B)$  to the category, $Vect$, of finite dimensional vector spaces over the field $\mathbb{K}$. However let us give here a more basic definition of a homotopy quantum field theory.

\subsection{Definition of HQFTs}
 A \emph{$(d + 1)$-dimensional homotopy quantum field theory, $\tau$, with background $B$}  assigns 
\begin{itemize}
\item to any $d$-dimensional $B$-manifold, $(X,g)$, a vector space, $\tau{(X,g)}$, 
\item to any $B$-isomorphism, $\phi : (X, g) \to ( Y, h)$, of $d$-dimensional $B$-manifolds, a $K$-linear isomorphism $\tau(\phi)  : \tau{(X, g)} \to \tau{( Y, h)}$,\\ 
\hspace*{-1cm} and
\item to any $B$-cobordism, $(W,F) : (X_0,g_0) \to (X_1,g_1)$, a $K$-linear transformation, $\tau(W,F) : \tau{(X_0,g_0)} \to \tau{(X_1,g_1)}$.
\end{itemize}
These assignments are to satisfy the following axioms:
\begin{enumerate}
\item[(1)] $\tau$ is functorial in $\mathbf{Man}(d,B)$, i.e., for two $B$-isomorphisms, $\psi: (X, g) \to ( Y, h)$ and $\phi : ( Y, h) \to (P,j)$, we have
$$\tau(\phi\psi) = \tau(\phi)\tau(\psi),$$
and if $1_{(X,g)}$ is the identity $B$-isomorphism on $(X,g)$, then $\tau(1_{(X,g)}) = 1_{\tau{(X,g)}}$.
\item[(2)]  There are natural isomorphisms
$$c_{(X,g),(Y,h)} : \tau((X,g)\amalg (Y,h)) \cong \tau(X,g)\otimes \tau(Y,h),$$
and an isomorphism, $u : \tau(\emptyset) \cong K$, that satisfy the usual axioms for a symmetric monoidal functor.
\item[(3)] For $B$-cobordisms, $(W,F) : (X,g) \to (Y,h)$ and $(V,G): (Y^\prime, h^\prime) \to (P,j)$ glued along a $B$-isomorphism $\psi :(Y,h) \to (Y^\prime,h^\prime)$, we have 
$$\tau((W,F)\amalg_\psi (V,G))= \tau(V,G)\tau(\psi)\tau(W,F).$$
\item[(4)] For the identity $B$-cobordism, $1_{(X,g)} = (I\times X, 1_g)$, we have 
$$\tau( 1_{(X,g)}) = 1_{\tau(X,g)}.$$
\item[(5)] For $B$-cobordisms $(W,F) : (X,g) \to (Y,h)$ and $(V,G) : (X^\prime,g^\prime) \to (Y^\prime,h^\prime)$ and $(P,J): \emptyset \to \emptyset$, the following diagrams are commutative:
$$\xymatrix{\tau((X,g)\amalg (X^\prime,g^\prime))\ar[r]^c\ar[d]_{\tau((W,F) \amalg (V,G))} &\tau(X,g)\otimes \tau(X^\prime,g^\prime)\ar[d]^{\tau(W,F) \otimes
\tau(V,G)} & &&\tau\emptyset \ar[r]^u\ar[d]_{\tau(P,J)}&K\\
\tau((Y,h)\amalg (Y^\prime,h^\prime))\ar[r]^c &\tau(Y,h)\otimes \tau(Y^\prime,h^\prime)&&&\tau\emptyset\ar[ur]_u&}.$$
\end{enumerate}

\textsc{Remark.}
These axioms are slightly different from those given in the original paper, \cite{turaev:hqft1}.  The really significant difference is in axiom 4 which is weaker than as originally formulated, where any $B$-cobordism structure on $I \times X$ was considered as trivial.  The effect of this change is important for us in as much as it is now the case that the HQFT is determined by the $(d+1)$-type of $B$, cf. Rodrigues, \cite{rodrigues}.  Because of this, it is feasible to attempt  a full classification of all $(1 + 1)$-HQFTs as there are simple algebraic models for 2-types, namely crossed modules.  We embarked on such a classification in \cite{TPVT:hqft1}. In general there is a challenge to relate algebraic models for the $(d+1)$-type of $B$ to the structure of the corresponding $(d+1)$-HQFTs, that is, a particular sort of classification problem.  Different algebraic models may give different perspectives on the HQFTs and a knowledge of the overall structure common to all the HQFTs with given background should allow us to identify any special features of particular HQFTs and thus to understand the geometry they are encoding.   

To be able to discuss classification of HQFTs, it is first necessary to discuss some notion of map between different such theories.

\medskip

\textsc{Definition.}
Let $\tau$ and $\rho$ be two $(d + 1)$-HQFTs with background $B$, then a map $\theta : \tau \to \rho$ is a family of maps $\theta(X,g) : \tau(X,g) \to \rho(X,g)$ indexed by the $B$-manifolds $(X,g)$ such that for every $B$-isomorphism, $\psi : (X,g) \to (Y,h)$, and every $B$-cobordism, $(W,F) : (X,g) \to (Y,h)$, the maps $\theta(X,g)$ and $\theta(Y,h)$ satisfy the obvious naturality conditions and conditions for compatibility with the structure maps, $r$, $l$, etc., (cf. Turaev's \cite{turaev:hqft1}, section 1.2, or Rodriques, \cite{rodrigues}, definition 1.4).

\medskip

Using this we can define a category $\mathbf{HQFT}(d,B)$ with obvious objects and maps.  Change of background space induces a functor between the corresponding categories and extending a result of Turaev (for the initial form of HQFT), Rodrigues proved in \cite{rodrigues} that the equivalence class of 
$\mathbf{HQFT}(d,B)$ depended only on the homotopy $(d+1)$-type of $B$.  One form of the classification problem mentioned above is thus to start with an algebraic model of the  $(d+1)$-type of $B$ and to find an algebraic description of the category $\mathbf{HQFT}(d,B)$ in terms of that algebraic structure.  For instance, if $B$ is  a $K(G,1)$, then Turaev proved that there is a bijective correspondence between the isomorphism classes of $(1 + 1)$-dimensional HQFTs with background $K(G,1)$ and isomorphism classes of crossed $G$-algebras (see \cite{turaev:hqft1}).  Brightwell and Turner, \cite{B&T}, for $B$ a $K(G,2)$ with, of course, $G$ Abelian, showed that $(1+1)$-dimensional HQFTs with such a background form a category equivalent to that of $G$-Frobenius algebras, i.e., Frobenius algebras with a specified $G$-action.  One of the aims of the first paper in this series, \cite{TPVT:hqft1}, was  to introduce algebraic objects that generalise both the crossed $G$-algebras and the $G$-Frobenius algebras and such that the categories of these objects would correspond to categories $\mathbf{HQFT}(1,B)$.  

\medskip

We note that a consequence of the definition of a homotopy quantum field theory is that if $\tau$ is a $(d+1)$-HQFT and $(X,g)$ and $(X,h)$ are two $B$-manifolds with the same underlying manifold, $X$, and the two characteristic maps, $g$ and $h$, are freely homotopic, then a choice of homotopy $F: g \simeq h$ gives a $B$-cobordism $(I\times X,F)$ which induces an isomorphism between $\tau(X,g) $ and $\tau(X,h)$. (This is an easy exercise, but is also a consequence of Rodrigues, \cite{rodrigues}, Proposition 1.2.) Because of this one should expect that some of the essential features of $\tau(X,g)$ will be influenced by the homotopy class of $g$.

 \section{Crossed gadgetry.}
\subsection{Crossed modules and their relatives.}
In the construction of examples of topological and homotopical quantum field theories, one often uses a (finite) group $G$, and a triangulation of the manifolds, $\Sigma$, etc., involved, and one assigns labels from $G$ to each (oriented)  edge of each (oriented) triangle, for example, (see diagram below),
\begin{figure}[h]\begin{center}
\font\thinlinefont=cmr5
\begingroup\makeatletter\ifx\SetFigFont\undefined
\def\x#1#2#3#4#5#6#7\relax{\def\x{#1#2#3#4#5#6}}%
\expandafter\x\fmtname xxxxxx\relax \def\y{splain}%
\ifx\x\y   
\gdef\SetFigFont#1#2#3{%
  \ifnum #1<17\tiny\else \ifnum #1<20\small\else
  \ifnum #1<24\normalsize\else \ifnum #1<29\large\else
  \ifnum #1<34\Large\else \ifnum #1<41\LARGE\else
     \huge\fi\fi\fi\fi\fi\fi
  \csname #3\endcsname}%
\else
\gdef\SetFigFont#1#2#3{\begingroup
  \count@#1\relax \ifnum 25<\count@\count@25\fi
  \def\x{\endgroup\@setsize\SetFigFont{#2pt}}%
  \expandafter\x
    \csname \romannumeral\the\count@ pt\expandafter\endcsname
    \csname @\romannumeral\the\count@ pt\endcsname
  \csname #3\endcsname}%
\fi
\fi\endgroup
\mbox{\beginpicture
\setcoordinatesystem units <0.750000cm,0.750000cm>
\unitlength=0.750000cm
\linethickness=1pt
\setplotsymbol ({\makebox(0,0)[l]{\tencirc\symbol{'160}}})
\setshadesymbol ({\thinlinefont .})
\setlinear
%
%
\linethickness= 0.500pt
\setplotsymbol ({\thinlinefont .})
\plot  3.6 21.220  5.5 23.732 /
%
%
\plot  5.394 23.531  5.5 23.732  5.309 23.594 /
%
%
%
\linethickness= 0.500pt
\setplotsymbol ({\thinlinefont .})
\plot  5.50 23.732  7.436 21.220 /
%
%
\plot  7.265 21.355  7.436 21.220  7.349 21.420 /
%
%
%
\linethickness= 0.500pt
\setplotsymbol ({\thinlinefont .})
\putrule from  3.573 21.139 to  7.436 21.139
%
%
\plot  7.224 21.086  7.436 21.139  7.224 21.192 /
%
%
%
\put{$g$} [lB] at  3.863 22.3
%
%
\put{$h$} [lB] at  6.8 22.3
%
%
\put{$k$} [lB] at  5.4 20.5
%
%
\put{ $\circlearrowleft$} [lB] at  5.15  22
%
%
\linethickness= 5.500pt
\setplotsymbol ({\thinlinefont .})
\ellipticalarc axes ratio  0.093:0.093  360 degrees 
	from  3.613 21.114 center at  3.520 21.114
\linethickness=0pt
\putrectangle corners at  3.410 23.757 and  7.461 20.701
\endpicture}

\label{figure 1}
\end{center}\end{figure}
with the boundary / cocycle condition that $kh^{-1}g^{-1} = 1$, so $k = gh$.

The geometric intuition behind this is that `integrating' a $G$-valued function around the triangle  yields the identity.  This intuition corresponds to problems where a $G$-bundle on $\Sigma$ is specified by charts and the elements $g$, $h$, $k$, etc. are transition automorphisms of the fibre.  The methods then use manipulations of the pictures as the triangulation is changed by subdivision, etc.

Another closely related view of this is to consider continuous functions $f : \Sigma \to BG$ to the classifying space of $G$.  If we triangulate $\Sigma$, we can assume that $f$ is a cellular map using a suitable cellular model of $BG$, at the cost of replacing $f$ by a homotopic map and perhaps subdividing the triangulation.  From this perspective the previous model is a \emph{combinatorial} model of such a continuous `characteristic' map, $f$.  The edges of the triangulation pick up group elements since the end points of each edge get mapped to the base point of $BG$, and $\pi_1BG \cong G$, whilst the faces give a realisation of the cocycle condition.
 
The natural first generalisation of this imagines a value assigned to the triangle itself, which measures the extent to which the cocycle condition is not satisfied. This  uses the concept of crossed module, which we recall next.
\medskip

\textsc{Definition.}
A \emph{crossed module}, $\mathcal{C} = (C,P, \partial)$, consists of groups $C$, $P$, a (left) action of
$P$ on $C$ (written $(p,c)\rightarrow {}^pc)$ and a homomorphism
$$\partial : C \rightarrow P$$
 such that

CM1 \quad $\partial({}^pc) = p\cdot\partial c\cdot p^{-1}$ \quad\quad for all $p 
\in P$, $c \in C$;\\
and

 CM2 \quad ~ ${}^{\partial c} c^\prime = c\cdot c^\prime\cdot c^{-1}$ \quad\quad\quad for all
 $c, c^\prime \in C$.

\medskip

Morphisms of crossed modules are pairs of maps preserving structure.  These give a category,
$\bf{CMod}$, of crossed modules.

We have here the definition for a crossed module of groups. It is fairly simple to adapt it to crossed modules of groupoids.  The only real points to note are that (i)  the two groupoids $C$ and $P$ have the same set of objects, (ii)  the morphisms are the identity on objects, and consequently, (iii) in $\mathcal{C} = (C,P, \partial)$ with $P$ and $C$  groupoids, $C$ is a `family of groups', that is,  if $x, y$ are distinct objects of $P$ and thus of $C$ as well, then $C(x,y)$, the set of arrows from $x$ to $y$ in $C$, will be empty.  The reason for this is clear when the motivating examples are considered as $C$ generalises both the idea of a normal subgroupoid and also the relative homotopy groups for varying base point. In the latter example, $C$ is the family $\{\pi_2(X,A,x)\}$ for $x$ in a set of basepoints and the action is the usual one coming from `change of base point' along a path.

\subsection{Internal Categories in the category of Group(oid)s.}

Let $\mathbf{C}$ be a category with finite limits, for instance, that of groups or Lie algebras.  An internal category in
$\mathbf{C}$ is a diagram\\
$$\xymatrix{C_1\ar[r]<1ex>^s\ar[r]<-1ex>_t&C_0\ar[r]^i &C_1 }.$$
and $si = ti = Id_{C_0}$,
together with a \emph{composition map} in $\mathbf{C}$,
$$C_1 ~_s \times_t C_1 \rightarrow C_1,$$
whose domain is given by the pullback $$\xymatrix{C_1 ~_s \times_t
  C_1\ar[r]\ar[d] &C_1\ar[d]^t\\C_1\ar[r]^s &C_0},$$
satisfying the usual associativity and identity rules.  We say $C_1$ is the
object of arrows and $C_0$ the object of objects, and then $C_1 ~_s \times_t C_1$
is the object of composable pairs of arrows.

We write $\mathbf{Cat(C)}$ for the category of internal categories in $\mathbf{C}$ and note:\\
\centerline{\emph{The categories $\mathbf{CMod}$ and $\mathbf{Cat(Grps)}$ are equivalent.}}
The proof, for instance, in \cite{RB&CBS}, is a generalisation of the very basic result from group theory that a congruence on a group is specified exactly by the normal subgroup of those elements equivalent to the identity element, whilst, conversely, any normal subgroup, $N \vartriangleleft G$ defines a congruence on $G$ by $g_1\sim g_2 $ if and only if $g_1g_2^{-1}$ is in $N$. The result seems to have been noticed by Grothendieck and  Verdier in the 1960s, but was not published until the paper \cite{RB&CBS} by Brown and Spencer.

There are other useful variant descriptions or interpretations of these same basic structures or of their lax analogues: (strict) 2-groups, categorical groups, cat$^1$-groups  and gr-categories being the most common ones.

\subsection{Crossed complexes}

The notion of a crossed complex of groups was defined by Blakers in 1946 (under the name of `group system') and  Whitehead (1949) considered free such objects that he called `homotopy systems'.  Blakers had used them as a way of systematising the properties of relative homotopy groups of a filtered space (and this may be an important perspective if it should prove useful to consider `extended HQFTs). Their theory has been developed extensively by Brown and Higgins in a long series of papers. They have also been considered by Baues, \cite{HJB:AH,HJB:4D}, playing a key role in his theory of combinatorial homotopy. In his work they are called \emph{crossed chain complexes}.

Let $$\mathbf{X}_\ast : X_0 \subseteq X_1\subseteq X_2\subseteq \cdots \subseteq X_n \subseteq \cdots \subseteq X_\infty$$
be a filtered space, then there are relative homotopy groups, $\pi_n(X_n,X_{n-1},x)$ for $x \in X_0$, obtained as relative homotopy classes of mappings of an $n$-cube into $X_n$ in such a way that all but one face of the $n$-cube goes to $x$, and the last face goes into $X_{n-1}$.  For $n > 1$, this gives a family of groups indexed by the basepoints $x \in X_0$.  For $n = 1$, the same idea yields the fundamental groupoid, 
$\pi_1(X_1X_0)$, of $X_1$ relative to $X_0$, so here the elements of the structure are end-point-fixed homotopy classes of paths in $X_1$ between points in $X_0$.  This structure at the lowest level of the filtration acts in a fairly obvious way (change of base point) on all the structures at higher level.  That is not the only structure linking these objects.  An element of $\pi_n(X_n,X_{n-1},x)$ is given as a relative homotopy class of mappings from $I^n$ into $X_n$ as mentioned above.  Restricting that mapping to the `last face' gives a mapping from $I^{n-1}$ to $X_{n-1}$ and this induced a homomorphism
$$\pi_n(X_n,X_{n-1},x)\stackrel{\partial}{\to}\pi_{n-1}(X_{n-1},X_{n-2},x),$$
which is compatible with the action of $\pi_1(X_1X_0)$.  This rich structure starts to look quite complex, but it can be abstracted and encoded moderately easily as follows:

\medskip

\textsc{Definition.}
A \emph{crossed complex} $\mathcal{C}$ (of groupoids) is a sequence of morphisms of groupoids over $C_0$
$$\xymatrix{\cdots \ar[r]&C_n\ar[r]^{\delta_n}\ar[d]^\beta &C_{n-1}\ar[r]^{\delta_{n-1}}\ar[d]^\beta &\cdots\ar[r]&C_2\ar[r]^{\delta_2}\ar[d]^\beta &C_1\ar[d]<.7ex>^{\delta^1}\ar[d]<-.7ex>_{\delta^0}\\
&C_0 &C_0&&C_0&C_0 }.$$
Here for each $n\geq 2$, $\{C_n\}$ is a family of groups with indexing  map $\beta$, so that for $p\in C_0$, we have groups $C_n(p) = \beta^{-1}(p)$. For $n=1$, $\delta^0$ and $\delta^1$ are the source and target maps of the groupoid $C_1$. Further we have an action of $C_1$ on each $C_n$, so that if $a:p\to q$ in $C_1$ and $x \in C_n(p)$, then ${}^ax\in C_n(q)$.  (This is a left action. Many sources use right actions, but it is more convenient for our situation to use a left action. Of course, the translation from one to the other is easy, if sometimes slightly confusing.) This data is to satisfy :

\renewcommand{\theenumi}{\roman{enumi}}
\begin{enumerate}
\item each $\delta_n$ is a morphism over the identity on $C_0$;
\item $\delta_2 : C_2\to C_1$ is a crossed module (over $C_1$);
\item $C_n$ is a $C_1$-module for $n \geq 3$;
\item $\delta_n : C_n \to C_{n-1}$ is a morphism compatible with the $C_1$-actions, for $n \geq 3$;
\item $\delta\delta : C_n\to C_{n-2}$ is trivial for $n \geq 3$; 
\item $\delta_2 C_2$ acts trivially on $C_n$ for $n \geq 3$.
\end{enumerate}

\textsc{Remarks.}
(i)  The family of groups $C_n$ can be also thought of as a groupoid with $C_0$ as its set of objects and if $p\in C_0$, $C_n(p,p) = C_n(p) = \beta^{-1}(p)$.  If $q \neq p\in C_0$, then $C_n(p,q) = \emptyset$. We will often be dealing with a \emph{reduced} crossed complex in which case $C_0$ is a singleton set and $C_n$ is just a single group.

(ii) If $C$ is a crossed complex, its fundamental groupoid $\pi_1C$ is the quotient of the groupoid, $C_1$, by the normal, totally disconnected subgroupoid $\delta_2C_2$.  By axiom (vi), the $C_n$ for $n \geq 3$ inherit a   $\pi_1C$-module structure and, as $\delta_2 : C_2\to C_1$ is a crossed module, there is an induced  $\pi_1C$-module structure on $\ker \delta_2$ as well.  This means that the $C_n$s, $\ker \delta_2$, and the boundary maps $\delta_n$, $n\geq 3$ yield a chain complex of  $\pi_1C$-modules.  This is most transparent when $C$ is reduced as then $\pi_1C$ is simply a group, and in dimensions greater than 2, we have the usual notion of a chain complex of modules over that group.

(iii) If $\mathbf{X}$ is a filtered space, we will denote by $\pi(\mathbf{X})$, the corresponding crossed complex.

\medskip

A \emph{morphism} $f : C \to D$ of crossed complexes is a family of groupoid morphisms, $f_n :C_n\to D_n$, for $n\geq 0$ compatible with boundaries and actions.  For $n = 0$, $f_0 $ is, of course, simply a function mapping the objects set $C_0$ to $D_0$ and all the groupoid morphisms, $f_n$, have $f_0$ as their object mapping.

This defines the category  $\mathbf{Crs}$ of crossed complexes.  This category has a very rich structure.  It is monoidal closed with an internal function-space functor, $\mathcal{CRS}$, and an associated tensor product. It also has an enrichment over the category of simplicial sets, and a closed model category structure in the sense of Quillen.  There is a notion of \emph{free crossed complex} and, in particular, if $X$ is a CW-complex and we give it the natural filtration by its skeleta, yielding a filtered space $\mathbf{X}$, then the construction sketched above will give us a free crossed complex, denoted $\pi(\mathbf{X})$. 

Triangulations of manifolds and other spaces or the \v{C}ech nerve construction lead to simplicial complexes and thus to simplicial sets if one puts an order on the vertices. Any simplicial complex or simplicial set, $K$,  has a geometric realisation, $|K|$. This space is naturally filtered by skeleta and we can therefore compose this with $\pi(-)$ to get a crossed complex. This assignment gives us a functor that also will be denoted $\pi $,  but this time $\pi : \mathcal{S} \to \mathbf{Crs}$, where $\mathcal{S}$ is the category of simplicial sets.  This functor $\pi : \mathcal{S} \to \mathbf{Crs}$ is left adjoint to a `nerve' functor, just as the geometric realisation is left adjoint to the singular complex functor.

\subsection{Nerves and classifying spaces of crossed complexes}\footnote{For readers not that accustomed to the setting of simplicial sets, we suggest a `classical' source namely the early chapters of Curtis' survey, \cite{curtis}. Other introductory material can be found in Kamps and Porter, \cite{K&P}, and there are also several  text books that describe other more recent developments, (e.g. Goerss and Jardine's \cite{GoerssJard}.)} (We assume the basic definitions of simplicial sets, simplicial groups, etc., are known.)

 For a crossed complex $\mathcal{C}$,
the simplicial set, $Ner(\mathcal{C})$, has a neat description as a `singular complex'.  

Let $\pi(n) = \pi(\Delta[n])$ be the free crossed complex on the $n$-simplex, $\Delta[n]$, then $Ner(\mathcal{C})_n \cong \mathbf{Crs}(\pi(n), \mathcal{C})$, the set of crossed complex maps from $\pi(n)$ to $\mathcal{C}$.  
\medskip

We will look in some more detail at the structure of $Ner(\mathcal{C})$, when $\mathcal{C}$ is, for simplicity,  a reduced crossed complex. (For the general case, we refer the reader to the discussion in \cite{ashley}, p.37, and \cite{NanTie1,NanTie2}.)   In the reduced case,  each of the families $C_n$ consists of a single group only, so \begin{itemize}
\item $Ner(\mathcal{C})_0$ consists of the set of  objects of $\mathcal{C}$, so is a singleton;
\item $Ner(\mathcal{C})_1 \cong C_1$, with, of course, $C_1$ considered merely as a set;
\item $Ner(\mathcal{C})_2 \cong (C_2\times C_1)\times C_1,$  and it thus has as typical element $(h_1,h_0)$, where $h_1 = (c_2,c_1)$ and $h_0 = c_1^\prime \in C_1$;
\item $Ner(\mathcal{C})_3 \cong (C_3\times C_2\times C_2\times C_1)\times (C_2\times C_1)\times C_1$ and has typical element, $(h_2,h_1,h_0)$, each $h_i$ being in the corresponding set as bracketted.
\end{itemize}
This way of viewing the simplicial set \emph{is} feasible in low dimensions and explicit formulae \emph{can} be derived, but is not that practical for $n\geq 4$.  In such dimensions it is often  more convenient to use a combination of the description as a `singular complex',  $Ner(\mathcal{C})_n \cong \mathbf{Crs}(\pi(n), \mathcal{C})$, together with an analysis, (cf.,   \cite{B&H1991}, p.99) of the structure of $\pi(n)$.  Each `singular simplex', $\sigma : \pi(n) \to \mathcal{C}$, will clearly select an element of $C_n$, together with a family of faces of decreasing dimensions.  For instance there will be a list of $n$ elements from $C_1$, corresponding to the maximum length path in $\sigma$ from the zeroth to the last vertex;  each adjacent pair in  these will determine a 2-simplex and so will need an element of $C_2$, to fill that simplex.  This will determine the third face of that 2-simplex.  Each triple from the one dimensional list, together with 2-dimensional fillers will give a partial shell of a 3-face of $\sigma$ and there will be a label from $C_3$ on the corresponding 3-face, and so on.  This iterative description was developed by Ashley, \cite{ashley}.

In general, $Ner(\mathcal{C})$ is a Kan complex, so any horn has a filler. In fact it has a stronger structure, namely that of a $T$-complex, again see \cite{ashley} and \cite{NanTie1,NanTie2} for a more detailed discussion of this.  The $T$ stands for `thin', cf. \cite{B&H1991}. In each dimension $n \geq 1$, there are certain elements called `thin'.  General elements of $Ner(\mathcal{C})$ can be specified by a morphism of crossed complexes $f : \pi(n) \to \mathcal{C}$.  Such an element will be called \emph{thin} if $f$ maps the top dimensional generating element of $\pi(n)$ to the/an identity element of $\mathcal{C}$.  These thin elements satisfy Dakin's axioms: \emph{degenerate elements are thin; any horn has a unique thin filler; if all faces but one of a thin element are thin, then so also is the last face}, (again see \cite{ashley}, \cite{B&H1991}, p.100, and references there).  Intuitively speaking the thin elements provide `canonical' fillers for horns.  There may be many different fillers for any horn (see the discussion on p.708 of \cite{tqft1}), but there will always be a unique \emph{thin filler}.
\footnote{A detailed study of the relationships between the $\overline{W}$ functor, $T$-complexes, and the nerve can be found in Nan Tie, \cite{NanTie1,NanTie2}.}  These $T$-complexes are particularly relevant in the context of this conference as they are closely related to the complicial sets introduced by Ross Street to handle coherence problems and weak $n$-categories in a closely related context, see the talk by Verity at this conference.

Any crossed complex $\mathcal{C}$ has a \emph{classifying space} $B\mathcal{C}$.  This can be obtained by taking the geometric realisation of $Ner(\mathcal{C})$, see \cite{B&H1991}. 

\subsection{Crossed Complexes and Simplicially Enriched Groupoids}
In order to justify that the general notion of formal map, to be introduced shortly, is a good combinatorial model for the characteristic maps needed for HQFT theory, it will be useful to recall briefly some of the basic links between simplicial objects and crossed complexes. Given any crossed complex $\mathcal{C}$, there is a simplicially enriched groupoid which can be constructed from it, and such that the corresponding Moore complex is isomorphic to $\mathcal{C}$.  This relationship between crossed complexes and `simplicial groupoids' will later give us important information on the geometric interpretation of our `extra structure' in this formal model, so we give a brief introduction to it here.

We denote the category of simplicial sets by $\mathcal{S}$. This is Cartesian closed so can be used as a base category in the sense of enriched category theory. In particular, Dwyer and Kan, and also Joyal and Tierney, introduced simplicially enriched groupoids which are many object analogues of simplicial groups. The category of such objects will be denoted by $\mathcal{S}\!-\!\mathbf{Gpds}$.\footnote{ We will often abbreviate the term `simplicially enriched groupoid' to $\mathcal{S}$-groupoid, but the reader should note that in some of the sources on this material the looser term `simplicial groupoid' is used to describe these objects usually with a note to the effect that this is not a completely accurate term to use.}

The loop groupoid functor of Dwyer and Kan, \cite{D&K},\label{loop-gpd} is a functor $$G: \mathcal{S} \to \mathcal{S}\!-\!\mathbf{Gpds},$$
which takes the simplicial set $K$ to the simplicially enriched groupoid $GK$, where $(GK)_n$ is the free groupoid on the directed graph
$$\xymatrix{K_{n+1}\ar[r]<.5ex>\ar[r]<-.5ex>&K_0},$$
where the two functions, $s$, source, and $t$, target, are $s = (d_1)^{n+1}$ and $t = d_0(d_2)^n$ with relations $s_0x = id$ for $x \in K_n$.  The face and degeneracy maps are given on generators by $$s_i^{GK}(x) = s_{i+1}^K(x),$$
 $$d_i^{GK}(x) = d_{i+1}^K(x),$$
for $x \in K_{n+1}$, $1 < i\leq n$ and $d_0^{GK}(x) = (d_1^K(x))(d_0^K(x))^{-1}.$ 

This loop groupoid functor has a right adjoint, $\overline{W}$, called the \emph{classifying space} functor. It generalises the nerve functor defined on small categories. This can be used to show that simplicially enriched groupoids model \emph{all} homotopy types, \cite{D&K}, extending the classical result that simplicial groups model all connected homotopy types. As for homotopy quantum field theory, it is natural to handle disjoint unions of manifolds, this is an important benefit for us as it implies that $G$ preserves disjoint unions and other colimits.

Given any $\mathcal{S}$-groupoid, $G$, its Moore complex $NG$ is given by 
$$NG_n = \bigcap_{i = 1}^n Ker(d_i : G_n \to G_{n-1})$$
with differential $\partial : NG_n \to NG_{n-1}$ being the restriction of $d_0$.  If $n \geq 1$, this is just a disjoint union of groups, one for each object in the object set, $O$, of $G$.  If we write $G\{x\}$ for the simplicial group of elements that start and end at $x \in O$, then at object $x$, one has
$$NG\{x\}_n = (NG_n)\{x\}.$$
In dimension 0, one has $NG_0 = G_0$, so the $NG_n\{x\}$, for different objects $x$, are linked by the actions of the 0-simplices, acting by conjugation via repeated degeneracies.

For simplicity in the description below, we will often assume that the $\mathcal{S}$-groupoid is \emph{reduced}, that is, its set $O$, of objects is just a singleton set $\{\ast\}$, so $G$ is just a simplicial group. We have used similar terminology  for crossed modules, crossed complexes, etc.

Suppose that $NG_m$ is trivial for $m >n$.

If $n =0$, then $NG_0$ is just the group $G_0$  and the simplicial group (or groupoid) represents an Eilenberg-MacLane space, $K(G_0,1)$.

If $n = 1$, then $\partial :NG_1 \to NG_0$ has a natural crossed module structure.

Returning to the discussion of the Moore complex, if $n = 2$,  then 
$$NG_2\stackrel{\partial}{\to}NG_1\stackrel{\partial}{\to}NG_0$$
has a 2-crossed module structure in the sense of Conduch\'e, \cite{conduche}. 

In all cases, the simplicial group will have homotopy groups only in the range covered by the non-trivial part of the Moore complex.  Using the Conduch\'e decomposition lemma, \cite{conduche}, one can decompose any level of any simplicial groupoid, $G$, as a semidirect product of factors given by terms of its Moore complex generalising the constructions used in the Dold-Kan theorem.  There is a reconstruction technique giving $G$, up to isomorphism, from $NG$ plus extra structure, a \emph{hypercrossed complex structure} in the sense of Carrasco and Cegarra, \cite{carrasco&cegarra}, see also the talk by Bourn at the conference. 

Now relaxing the restriction on $G$, for each $n > 1$, let $D_n$ denote the subgroupoid of $G_n$ generated by the degenerate elements. When handling $n$-types, for finite $n$, one can ask that Moore complex terms be trivial above some level, but instead of
asking that $NG_n$ be trivial, we can ask that $NG_n \cap D_n$ be. The importance of this is that the structural information on the homotopy type represented by $G$ includes structure such as the Whitehead products and these all lie in the subgroupoids $NG_n \cap D_n$. If these are all trivial then the algebraic structure of the Moore complex is simpler, being that of a crossed complex, and $\overline{W}G$ is isomorphic to  the \emph{nerve} of the crossed complex so its geometric realisation is the \emph{classifying space of that crossed complex}, cf. \cite{B&H1991}.  
 The crossed complex associated to a simplicial groupoid, $G$, is given explicitly by
 $$C(G)_n = \frac{NG_{n-1}}{(NG_{n-1}\cap D_{n-1})d_0(NG_{n}\cap D_{n})} \textrm{ \quad for  } n \geq 2,$$ 
$C(G)_1 = NG_0$, and, of course, $C_0$ is the common set of objects of $G$.

 
\section{Formal Maps}
\subsection{Formal $\mathcal{C}$-maps in low dimensions.}
In \cite{TPVT:hqft1}, we introduced the notion of a formal $\mathcal{C}$-map on 1- and 2-dimensional manifolds and also  on simplicial complexes.  Here $\mathcal{C}$ is a crossed module, $\mathcal{C}= (C,G,\partial)$:

\textsc{Definition, }\cite{TPVT:hqft1}.
Let $K$ be a simplicial complex.  A simplicial formal $\mathcal{C}$-map, $\lambda$, on $K$ consists of families of elements\\
(i) $\{c_t\}$ of $C$, indexed by the set, $K_2$, of 2-simplices of $K$, 
\\(ii) $\{p_e\}$ of $P$, indexed by the set of 1-simplices, $K_1$, of $K$,\\
 and a partial order on the vertices of $K$, so that each simplex is totally ordered. (This replaces the orientation and gives start vertices to all edges and triangles without problem.)  The assignments of $c_t$ and $p_e$, etc. are to satisfy\\
(a) the boundary condition  $$\partial c_t = p_1 p_0^{-1}p_2^{-1},$$
where the vertices of $t$, labelled $v_0$, $v_1$, $v_2$ in order, determine the numbering of the opposite edges, e.g., $e_0$ is between $v_1$ and $v_2$, and $p_{e_i}$ is abbreviated to $p_i$;\\
and \\
(b) the cocycle condition:\\
in a tetrahedron
yielding two composite faces
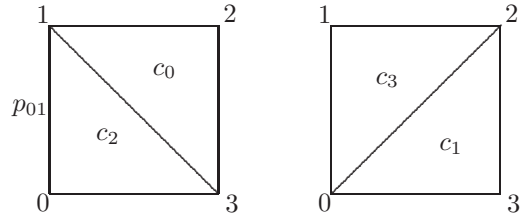
\begin{figure}[h]\begin{center}$\begin{array}{cc}
\font\thinlinefont=cmr5
\begingroup\makeatletter\ifx\SetFigFont\undefined
\def\x#1#2#3#4#5#6#7\relax{\def\x{#1#2#3#4#5#6}}%
\expandafter\x\fmtname xxxxxx\relax \def\y{splain}%
\ifx\x\y   
\gdef\SetFigFont#1#2#3{%
  \ifnum #1<17\tiny\else \ifnum #1<20\small\else
  \ifnum #1<24\normalsize\else \ifnum #1<29\large\else
  \ifnum #1<34\Large\else \ifnum #1<41\LARGE\else
     \huge\fi\fi\fi\fi\fi\fi
  \csname #3\endcsname}%
\else
\gdef\SetFigFont#1#2#3{\begingroup
  \count@#1\relax \ifnum 25<\count@\count@25\fi
  \def\x{\endgroup\@setsize\SetFigFont{#2pt}}%
  \expandafter\x
    \csname \romannumeral\the\count@ pt\expandafter\endcsname
    \csname @\romannumeral\the\count@ pt\endcsname
  \csname #3\endcsname}%
\fi
\fi\endgroup
\mbox{\beginpicture
\setcoordinatesystem units <0.5000cm,0.5000cm>
\unitlength=0.5000cm
\linethickness=1pt
\setplotsymbol ({\makebox(0,0)[l]{\tencirc\symbol{'160}}})
\setshadesymbol ({\thinlinefont .})
\setlinear
%
%
\put{2} [lB] at  9.129 22.754
%
%
\linethickness= 0.500pt
\setplotsymbol ({\thinlinefont .})
\plot  4.498 22.648  8.970 18.203 /
%
%
\put{0} [lB] at  4.155 17.727
%
%
\put{1} [lB] at  4.155 22.727
%
%
\linethickness= 0.500pt
\setplotsymbol ({\thinlinefont .})
\putrule from  4.498 18.176 to  4.498 22.648
\putrule from  4.498 22.648 to  8.996 22.648
\putrule from  8.996 22.648 to  8.996 18.176
\putrule from  8.996 18.176 to  4.498 18.176
%
%
\put{3} [lB] at  9.182 17.700
%
%
\put{$p_{01}$} [lB] at  3.520 20.373
%
%
\put{$c_2$} [lB] at  5.743 19.605
%
%
\put{$c_0$} [lB] at  7.250 21.351
\linethickness=0pt
\putrectangle corners at  3.520 22.754 and  9.390 17.700
\endpicture}
\quad&\quad
\font\thinlinefont=cmr5
\begingroup\makeatletter\ifx\SetFigFont\undefined
\def\x#1#2#3#4#5#6#7\relax{\def\x{#1#2#3#4#5#6}}%
\expandafter\x\fmtname xxxxxx\relax \def\y{splain}%
\ifx\x\y   
\gdef\SetFigFont#1#2#3{%
  \ifnum #1<17\tiny\else \ifnum #1<20\small\else
  \ifnum #1<24\normalsize\else \ifnum #1<29\large\else
  \ifnum #1<34\Large\else \ifnum #1<41\LARGE\else
     \huge\fi\fi\fi\fi\fi\fi
  \csname #3\endcsname}%
\else
\gdef\SetFigFont#1#2#3{\begingroup
  \count@#1\relax \ifnum 25<\count@\count@25\fi
  \def\x{\endgroup\@setsize\SetFigFont{#2pt}}%
  \expandafter\x
    \csname \romannumeral\the\count@ pt\expandafter\endcsname
    \csname @\romannumeral\the\count@ pt\endcsname
  \csname #3\endcsname}%
\fi
\fi\endgroup
\mbox{\beginpicture
\setcoordinatesystem units <0.5000cm,0.5000cm>
\unitlength=0.5000cm
\linethickness=1pt
\setplotsymbol ({\makebox(0,0)[l]{\tencirc\symbol{'160}}})
\setshadesymbol ({\thinlinefont .})
\setlinear
%
%
\linethickness= 0.500pt
\setplotsymbol ({\thinlinefont .})
\plot  4.551 18.212  9.023 22.657 /
%
%
\put{0} [lB] at  4.155 17.727
%
%
\put{1} [lB] at  4.155 22.727
%
%
\put{2} [lB] at  9.129 22.754
%
%
\put{3} [lB] at  9.182 17.700
%
%
\put{$c_3$} [lB] at  5.690 21.114
%
%
\put{$c_1$} [lB] at  7.383 19.340
%
%
\linethickness= 0.500pt
\setplotsymbol ({\thinlinefont .})
\putrule from  4.498 18.176 to  4.498 22.648
\putrule from  4.498 22.648 to  8.996 22.648
\putrule from  8.996 22.648 to  8.996 18.176
\putrule from  8.996 18.176 to  4.498 18.176
\linethickness=0pt
\putrectangle corners at  4.155 23.072 and  9.390 17.700
\endpicture}
\end{array}$
\end{center}\end{figure}\\
we have $$c_2{}^{p_{01}}c_0 = c_1c_3.$$

\subsection{Simplicial formal maps}
This notion is a `bare-hands' version of a more general one which can be given in terms of simplicial theory and which applies to a general crossed complex, $\mathcal{C}$, as base / coefficients. This latter approach is therefore more appropriate for modelling characteristic maps in higher dimension HQFTs and for developing general theory.

Let $K$ be a simplicial complex and $\mathcal{C} = (C_i, \partial_i)$, a reduced crossed complex. 

\medskip

\textsc{Definition.}
A \emph{simplicial formal $\mathcal{C}$-map} on $K$ is a pair consisting of an ordering, $\leq $, on the vertices of $K$, so that each simplex is totally ordered, and a simplicial map
$$\lambda : K \to Ner(\mathcal{C}).$$

\textsc{Remarks.}
(i) Note the ordering $\leq$ on $K_0$ endows $K$ with the structure of a simplicial set, which we will usually also write as $K$.  The simplices of the \emph{simplicial set} $K$ are the ordered sets, $\langle v_0\leq \ldots\leq  v_k\rangle$, where, after deletion of any repetitions, the resulting set, $\{v_0, \ldots, v_k\}$, is a simplex of the simplicial complex $K$, cf., for example, \cite{curtis} p.111. 

(ii)  The term `formal map' is suggested as recalling two images.  The first is that of a formally defined `mapping' from the realisation of $K$ to the classifying space, $B\mathcal{C}$, of $\mathcal{C}$, and we will explore this in more detail later. The second is that of a map on a surface, being an embedded graph with complement a disjoint union of discs, as in the idea of \emph{coloring} a map on a surface with elements of a group or other structure.

\medskip

Itemising the data specifying $\lambda$, using our earlier description of $Ner(\mathcal{C})$,  we find:\begin{itemize}
\item to each vertex of $K$, $\lambda$ assigns the single vertex of $Ner(\mathcal{C})$;
\item to each edge / 1-simplex $\sigma^{(1)}$ of $K$, an element $\lambda(\sigma^{(1)})$ of $X_1$, so, for instance, if $\mathcal{C}$ is a crossed module, as above, this will give an element of the group $G$;
\item to each 2-simplex, $\sigma = \sigma^{(2)}$, of $K$ an element of $C_2$, such that, abbreviating $d_i\sigma $ to $\sigma_i$,
$$\partial \lambda(\sigma) = \lambda(\sigma_1)\lambda(\sigma_2)^{-1} \lambda(\sigma_0)^{-1},$$
etc.
\end{itemize}
Thus we have a picture
$$\xymatrix{&\hspace{1mm}\ar[dr]^{\lambda(\sigma_0)}&\\
\hspace{1mm}\ar[ur]^{\lambda(\sigma_2)}_{\hspace{.3cm}\lambda(\sigma)}\ar[rr]_{\lambda(\sigma_1)}&&\hspace{1mm}}$$
and similarly into higher dimensions \footnote{The anticlockwise orientation in the formula, relative to the picture, is due to the various conventions, e.g. left actions, form of the Moore complex, etc. that we have made and should not be thought of as being important in itself.}. Restricting to $\mathcal{C}$ being a crossed module (so $C_n =1$ for $n\geq 3$), the cocycle condition comes from the fact that, if $\sigma = \sigma^{(3)}$ is a 3-simplex of $K$, then $\lambda(\sigma) = 1$.  This notion then reduces to our previous one.
 
(iii) $Ner$ is right adjoint to $\pi$, the functor from simplicial sets to crossed complexes, or if you prefer, to the Dwyer-Kan `simplicial groupoid' construction mentioned earlier,  so we could specify $\lambda$ by a map 
$$\overline{\lambda} : \pi(K) \to \mathcal{C},$$
or alternatively by a morphism of simplicially enriched groupoids 
$$\overline{\lambda} : GK \to \mathcal{C},$$
where we are thinking of the crossed complex, $\mathcal{C}$, as the corresponding simplicial group.

\medskip

Suppose that $\lambda : K \to Ner(\mathcal{C})$ is a formal map with an order $\leq $ given on the vertices of $K$. What happens if the order of the vertices is changed?  We will see later on that the assignment, $\lambda$, can be changed to be compatible with the new ordering. This will drop out of results on the way such a `coloring' behaves under subdivision. Clearly any formal map $\lambda : K \to Ner(\mathcal{C})$ induces a continuous mapping on the realisations, i.e., 
$$|\lambda| : |K| \to B\mathcal{C}$$and, in fact, the homotopy class of this is independent of the ordering on the vertices.

\subsection{Formal maps on a manifold.}
As we really want to look at manifolds, we can as a first step easily extend the idea of a formal map to one on a manifold relative to a triangulation $T$.

Let $X$ be a $d$-manifold, or more generally a $d$-dimensional polyhedron, and $\mathbf{T} = (T,\phi : |T|\to X)$ be an ordered triangulation of $X$, so $\phi $ is a homeomorphism between the realisation of the simplicial complex $T$ and $X$.  

\textsc{Definition.}
A \emph{(simplicial) formal  $\mathcal{C}$-map on $X$  relative to} $\mathbf{T}$ is a formal map $\lambda :T \to Ner(\mathcal{C})$ and hence, notationally, may be specified by $(T,\phi,\lambda)$ or, more briefly, $(\mathbf{T},\lambda)$.

It is sometimes important to remember that the data for a triangulation includes explicit mention of the homeomorphism, as the simplicial complex $T$ by itself is not enough to specify $\mathbf{T}$.  Alternatively, the simplicial complex may arise as the \v{C}ech nerve of an open cover of $X$ and the extra data needed is an ordering on the open sets making up the open cover.  This is sometimes useful for applications and also for links with cohomology, \cite{tqft1} and we will reurn to it in the last section of this paper.

Although our manifolds are oriented, the orientation only needs careful attention occasionally as many of the constructions do not use it explicitly, being special cases of more general ones.

\medskip

In addition to formal maps we will need formal cobordisms between them.  This can be done in more generality than we will give here, where we restrict to manifolds, as the notion of cobordism between manifolds is well known and well understood.  (A suitable setting for the extension to complexes can be given using, for instance, the domain categories of Quinn, \cite{quinn:atqft}.) We will need to consider triangulated manifolds $X_i,$ $ i = 1,2$, with triangulations $T_i$, and a cobordism $M$ between them, triangulated by $\mathcal{T}$ compatibly with the incoming and outgoing boundary triangulations. Of course, this implicitly involves the isomorphisms between these components and the two original manifolds, and the compatibility of the triangulating homeomorphisms with this structure.  We will not give this explicitly. 

 If $\lambda_1 : T_1 \to Ner(\mathcal{C})$ and $\lambda_2 : T_2 \to Ner(\mathcal{C})$ are 	two formal maps on the manifolds $X_1$ and $X_2$, then a \emph{formal $\mathcal{C}$-cobordism}, $\mathbf{\Lambda} : \lambda_1 \to \lambda_2$, consists of a triangulated cobordism, $(M,\mathcal{T})$, between them, and a formal $\mathcal{C}$-map, $\mathbf{\Lambda} : \mathcal{T}\to Ner(\mathcal{C})$, defined compatibly with the $\lambda_i$ on the incoming and outgoing boundaries. (Again we will not give this condition explicitly.)  We will usually be concerned with such formal cobordisms up to equivalence relative to the boundaries, in a sense that will be made precise shortly.

In the definition of simplicial formal map, we have taken the background to be $Ner(\mathcal{C})$ for $\mathcal{C}$ a reduced crossed complex. We could equally well have taken, as background here, any $BG = |\overline{W}G|$ for a general simplicial group $G$, and  this may be an important generalisation to make, however it is not completely clear how to generalise certain aspects of  the  formal maps when taking such a more general background and so we have made the restriction from the start.  In particular, using crossed complexes we can take account of any available cellular structure as follows.
\subsection{Cellular formal maps}
If we are considering a regular CW-decomposition of a space $X$, then there is an obvious generalisation of simplicial formal maps, extending the notion of cellular formal maps introduced in \cite{TPVT:hqft1}.

\textsc{Definition.}
 A \emph{(cellular) formal $\mathcal{C}$-map} on a  regular CW-complex $X$ is a crossed complex morphism 
$$\lambda : \pi(\mathbf{X}) \to \mathcal{C},$$
where, as before, $\mathbf{X}$ denotes the space $X$ with the skeletal filtration on it.
 
If the CW-structure came from a triangulation of $X$ then this coincides with the previous definition.
We will often omit the qualifying `simplicial' or `cellular', in the term `simplicial (or cellular) formal map'.  `Formal map' can therefore refer to either situation without real ambiguity.  One reason for working with crossed complexes rather than simplicial groups is that the transition from a CW-structure or any similar cell or handle decomposition of a space to the algebraic model (crossed complex) does not need the intervention of a supplementary triangulation, followed by elimination of the spurious effects that imposing that the triangulation brings.  We can usually go directly to the algebraic model of the geometry. 
\subsection{Equivalence of simplicial formal maps}
The idea behind the definition of a formal map is that it provides a good approximation to a characteristic map of a $B$-manifold, but is specified in an algebraic/combinatorial form.  It, of course, needs the triangulation or regular CW-decomposition on the underlying space, but clearly we will need to be able to subdivide or combine simplices or cells to get new decompositions and new `equivalent' formal maps. More precisely, suppose $\lambda : K \to Ner(\mathcal{C})$ is a formal map on an ordered triangulation of the space, $X$, and $K^\prime$ is another ordered triangulation of $X$, we need to have a notion of equivalent formal $\mathcal{C}$-maps and a technique for constructing such maps, so that (i) we can construct a new formal $\mathcal{C}$-map $\lambda^\prime$ on $K^\prime$ and that (ii) $\lambda$ and $\lambda^\prime$ are `equivalent'.

In \cite{TPVT:hqft1}, and thus for the case of 1-dimensional spaces, and `surfaces as cobordisms', but  with $\mathcal{C}$ being a crossed module, we achieved this by using the 3-dimensional cocycle condition, however in general that is not available to us so we need to use an alternative method.  We will mimic the treatment of \cite{TPVT:hqft1}, wherever possible and will sometimes just quote results if the proof goes across to the general case. This will allow us to indicate more clearly where the differences occur.

Although we have not yet defined formal HQFTs in this generality, we will put ourselves in the context that will be needed later by supposing that $X$ is a polyhedron with a given family of base points $\mathbf{m} = \{m_i\}$.  This will correspond either to having at least one basepoint in each connected component of the object or in each boundary component if $X$ is a cobordism between two objects.)  Let $K_0$, $K_1$ be two triangulations of $X$, i.e., $K_0$ and $K_1$ are simplicial complexes with geometric realisations homeomorphic to $X$ (by specified homeomorphisms) with the given base points among the vertices of the triangulation. 

\textsc{Definition.}
 Given two formal $\mathcal{C}$-maps $(K_0,\lambda_0),$ $(K_1,\lambda_1)$, then we say they are \emph{equivalent}  if there is a triangulation, $T$, of $X\times I$ extending $K_0$ and $K_1$ on $X \times \{0\}$ and $X\times \{1\}$ respectively, and a formal $\mathcal{C}$-map, $\Lambda$, on $T$ extending the given ones on the two ends and respecting the base points, in the sense that $T$ contains a subdivided $\{m_i\}\times I$ for each basepoint $m_i$ and $\Lambda $ assigns the identity element $1_P$ of $P$ to each 1-simplex of $\{m_i\}\times I$.

We will use the term `ordered simplicial complex' for a simplicial complex, $K$, together with a partial order on its set of vertices such that the vertices in any simplex of $K$ form a totally ordered set.  If we give the unit interval, $I$, the obvious structure of an ordered simplicial complex with $0 < 1$, then the cylinder $|K|\times I$ has a canonical triangulation as an orderd simplicial complex and we will write $K\times I$ for this.

If we are given two formal $\mathcal{C}$-maps defined on the same $K$, $(K,\lambda_0),$ and $(K,\lambda_1)$, we say they are \emph{simplicially homotopic} if there is a formal $\mathcal{C}$-map defined on the simplicial complex $K\times I$ extending them both.

In particular, when we are considering two formal $\mathcal{C}$-cobordisms $\mathbf{\Lambda}_i : \mathcal{T}_i\to Ner(\mathcal{C})$, $i = 1,2$, between two formal $\mathcal{C}$-maps, $\lambda_i : T_i \to Ner(\mathcal{C})$, (so $|\mathcal{T}_1| = |\mathcal{T}_2| $), we will need to consider \emph{(simplicial) homotopy relative to the boundaries} in the sense that the formal map on $|\mathcal{T}_1|\times I$, which gives the homotopy is constant on the two subcomplexes triangulating the ends of $|T_i|\times I$.
This is related to an obvious form of \emph{equivalence relative to the boundaries}.  In the following results the extension to the relative case is easy.
\begin{lemma} \label{refl-trans}
Equivalence is an equivalence relation.
\end{lemma}

\textsc{Proof.}
This is routine and extends without problem the proof in the low dimensional case of \cite{TPVT:hqft1}. Transitivity and symmetricity are easy, whilst reflexivity merely requires the construction of the obvious triangulation $\mathcal{K}$ of $X\times I$, followed by the obvious construction of a formal map on $\mathcal{K}$. The details are omitted. \hfill $\square$

Equivalence combines  the intuition of the geometry of triangulating a (topological) homotopy,  where the triangulations of the two ends may differ,  with some idea of a  combinatorially defined simplicial homotopy of formal maps. The proof of the following is immediate from the definition and is omitted.

\begin{lemma}If   $(K,\lambda_0)$ and $(K,\lambda_1)$ are two formal maps, which are simplicially homotopic as formal $\mathcal{C}$-maps, then they are equivalent. \hspace*{1cm}\hfill $\square$
\end{lemma}
This applies not only to the basic formal maps but to cobordisms between such maps.

\begin{proposition}\label{reorder}A change of order on the vertices of $K$  generates an equivalent formal $\mathcal{C}$-map.
\end{proposition}
\textsc{Proof.} 
Let $K_0$ be $K$ with the given order and $K_1$ the same simplicial complex with a new ordering.  Construct a triangulation $T$ of $|K|\times I$ having $K_0$ and $K_1$ on the two ends.  (Inductively, we can suppose just one pair of elements has been transposed in the order.)  Extend any given $\lambda_0$ on $K_0$ over $T$ and then restrict to get an equivalent $\lambda_1$ on $K_1$.\hfill $\square$

\begin{theorem}Given a simplicial complex, $K$, with geometric realisation $X = |K|$, and a subdivision $K^\prime$ of $K$.\\
(a)  Suppose $\lambda$ is a formal $\mathcal{C}$-map on $K$, then there is a formal $\mathcal{C}$-map, $\lambda^\prime$ on $K^\prime$ equivalent to $\lambda$.\\
(b) Suppose $\lambda^\prime$ is a formal $\mathcal{C}$-map on $K^\prime$, then there is a formal $\mathcal{C}$-map, $\lambda$ on $K$ equivalent to $\lambda^\prime$.
\end{theorem}

\textsc{Proof.}
 We triangulate a copy of the cylinder $|K|\times I$, so that we have the triangulation $K$ on $|K|\times \{0\}$ and $K^\prime$ on $|K|\times \{1\}$. (An explicit way of doing this is discussed in \cite{TPVT:hqft1}.) This can be done so that the simplices in $K$ that are unaffected by the subdivision yield prisms with the standard simplicial set structure.  In particular we have that the base points $m_i$ give 1-simplices $m_i \times I$ in the triangulated cylinder.

This set up is the same for both parts of the proof. Now assume given $\lambda $ defined on $K$ and thinking of $K$ as $K\times \{0\}$, we seek to extend $\lambda$ to a formal $\mathcal{C}$-map, say $\Lambda$, on the triangulated cylinder.  If we can do that we will be able to restrict $\Lambda$ to the copy of $K^\prime$ on $|K|\times \{1\}$ to get a formal $\mathcal{C}$-map, $\lambda^\prime$, and, by definition, this will be equivalent to $\lambda$ proving (a).  Reversing the roles of the two ends a similar argument will prove (b).

It thus remains to check that the extension $\Lambda$ exists.

Away from the extra vertices, we can extend $\lambda$ in an obvious way.  Use the value $1_{C_1}$ as the label for any vertical 1-simplex, and suppose that $\sigma$ is a $n$-simplex in both $K$ and $K^\prime$. The triangulation `above' $\sigma$ will be $\sigma \times \Delta[1]$ so we use the projection to $\sigma$ and the labelling $\lambda(\sigma)$ to define $\Lambda$ on $\sigma \times \Delta[1]$.  If, on the other hand, $\sigma$ in $K$ is subdivided in $K^\prime$, then we have the prism $|\sigma|\times I$ and the triangulation described in \cite{TPVT:hqft1}, consists of the joins of initial segments of $\sigma$ in the base with the complementary segment subdivided as necessary in the top. The extension scheme used in \cite{TPVT:hqft1} is easily adapted to our context by replacing the use of the cocycle condition by the Kan filling/extension condition.  The only other difference is that whilst in \cite{TPVT:hqft1} the extension process stopped in dimension 3, now we continue that extension process until  dimension $dim(K) +1$ is reached. (The reader is referred to \cite{TPVT:hqft1} for a more detailed description of this extension process in low dimensions.)

Finally this method does not depend on which end of the cylinder is used first so it is easily adapted to handle (b). \hspace*{1cm}\hfill $\square$

\textsc{Remarks.}
(i) As is usual with Kan complexes, we can think of filling simplices or extending maps as generalised or weak compositions. Thus using the Kan property of $Ner(\mathcal{C})$, we can compose values of a formal map on adjacent simplices.  As we have unique canonical `thin' fillers for all horns in $Ner(\mathcal{C})$, these compositions could in principle be written down exactly.  In fact some elementary cases of this process were given in \cite{TPVT:hqft1}.

(ii)  There is an alternative proof of the extension part of the above result.  It is very neat but less constructive so does not suggest that the composition process is algebraic as does the one used above: consider the diagram
$$\xymatrix{K\ar@{^{(}->}[d]\ar[r]^\lambda & Ner(\mathcal{C})\ar[d]\\
	         L \ar@{-->}[ur]^\Lambda\ar[r] & 1}$$
where $L$ denotes the triangulated cylinder.  The inclusion of the end $K$ into $L$ is a trivial cofibration, and $Ner(\mathcal{C})\to 1$ is a Kan fibration, so the dashed diagonal exists as required.

Given any cellular formal $\mathcal{C}$-map, we can triangulate the cell complex and find a simplicial formal $\mathcal{C}$-map that is cellularly equivalent to it.  Conversely given a simplicial formal map, $\lambda$, on a triangulation of a regular CW-complex, then we can `integrate' $\lambda$ over each cell, inductively up the skeleton, to get a cellular formal map equivalent to it.  The process in each case is to decompose the cylinder on the complex compatibly with the CW decomposition on one end and the triangulation on the other.  The extension argument given above does need adapting slightly as we are now in a more topological,  less simplicial, setting,  but the idea is essentially the same.  An example of this in low dimensions is given in \cite{TPVT:hqft1}.

\section{Formal maps as models for $B\mathcal{C}$-manifolds}
We next indicate why it is reasonable to expect the combinatorial mechanism of formal maps accurately to reflect the notion of a map from a polyhedral space or manifold to $B= B\mathcal{C}$, where, as before, $\mathcal{C}$ is a reduced crossed complex. 

\subsection{From `formal' to `actual'}
Given any formal $\mathcal{C}$-map
$$\lambda : K \to Ner(\mathcal{C}),$$ we can take its geometric realisation to get a map
$$|\lambda| : |K| \to |Ner(\mathcal{C})| = B\mathcal{C}.$$
We thus have a $B\mathcal{C}$-space and if, for instance, $K$ was an ordered triangulation of a manifold $M$, we could compose with the homeomorphism $\phi $, say, between $|K|$ and $M$ to get a $B\mathcal{C}$-manifold or cobordism.  It is clear that other choices of $\phi$ correspond to the action of the automorphism group of $M$ on the set of maps from $M$ to $B\mathcal{C}$  and so are already accounted for in the theory.

If $\lambda : K \to Ner(\mathcal{C})$ and $\lambda^\prime : K^\prime \to Ner(\mathcal{C})$ are equivalent formal maps, then the equivalence (i.e. the formal map on the cylinder) gives a reversible $B\mathcal{C}$-cobordism between the two resulting  $B\mathcal{C}$-manifolds. Again this is accounted for within the HQFT.

Going from `formal' maps to `actual' maps thus causes no problems. One just uses geometric realisation. To go in the other direction one expects to use simplicial and cellular approximation theory.


\subsection{Simplicial and CW-approximations and the passage to Crossed Complexes}

Suppose $K$ is an $n$-dimensional simplicial complex, then simplicial / CW-approximation theory implies that the space of maps from $|K|$ to $B\mathcal{C}$ is weakly homotopy equivalent to $|\mathcal{S}(K, Ner(\mathcal{C}))|$.
Thus any characteristic map $g : |K|\to B\mathcal{C}$ is in the same connected component  of this mapping space as a realisation, $|\lambda|$, of a formal $\mathcal{C}$-map.  Moreover any two ways of connecting $g$ to such a $|\lambda|$ will be mirrored by a pair of paths in $|\mathcal{S}(K, Ner(\mathcal{C}))|$.

 The simplicial set $\mathcal{S}(K, Ner(\mathcal{C}))$ is itself equivalent to $Ner(\mathcal{CRS}(\pi K, \mathcal{C}))$, where $\mathcal{CRS}(\mathcal{C}, \mathcal{D})$ denotes the Brown-Higgins crossed complex of morphisms from a crossed complex $\mathcal{C}$ to another one $\mathcal{D}$. This means that $|\mathcal{S}(K, Ner(\mathcal{C}))|$ is weakly equivalent to the classifying space of $\mathcal{CRS}(\pi K, \mathcal{C})$.  (These results are special cases of results of Brown and Higgins in the papers, \cite{B&H1987,B&H1991}.)

Now suppose that we are considering a $(d+1)$-HQFT, $\tau$, then all the objects, manifolds and cobordisms have dimension less than or equal to $d+1$.  As a consequence the images of all formal $\mathcal{C}$-maps will be trivial in dimensions greater than $d + 1$, since if $K$ has dimension $n$, the crossed complex, $\pi(\mathbf{K})$ will be trivial in dimensions greater than $n$, as $\pi(\mathbf{K})_p = \pi_p(K_p,K_{p-1},K_0)$ and is a free $\pi_1(K_1K_0)$-module on the $p$-cells of  $K$. We may, thus, replace $\mathcal{C}$ by its $(d+1)^{th}$ truncation, $tr_{d+1}\mathcal{C}$, without loss of generality.  To do this we replace each $C_n$ by the trivial group above dimension $d+1$ and replace $C_{d+1}$ by $C_{d+1}/\partial C_{d+2}$. Any formal $\mathcal{C}$-map on $K$ corresponds uniquely to  a formal $tr_{d+1}\mathcal{C}$-map and conversely. This is the analogue in this context of the fact that a general $d+1$-HQFT with background $B$ only depends, up to isomorphism, on the $(d+1)$-type of $B$.

The setting is now clear when it comes to equivalence of formal $\mathcal{C}$-cobordisms. If we have two equivalent formal $\mathcal{C}$-cobordisms between two formal $\mathcal{C}$-maps, then the equivalence corresponds to a $(d+2)$-dimensional simplicial complex in the form of a cylinder, (so the  highest dimensional simplices must be labelled by the identity elements of $C_{d+2}$ and hence correspond to a cocycle condition in this dimension).  As a result, the two induced maps under geometric realisation will be homotopic and the resulting induced maps under the HQFT will be equal. 


\textsc{Remarks.}
(i) Note that $\pi(\mathbf{K})$ can be given as a colimit, over the category of simplices of $K$, of the various $\pi(k)$, that is the crossed complex of a $k$-dimensional simplex.  For each $k$, and each $k$-simplex $\sigma \in K_k$, a formal map $\lambda$  yields a map from $\pi (k)$ to $\mathcal{C}$ and thus specifies an element in $\mathcal{C}_k$.  These different elements are related by face formulae to the corresponding elements in $\mathcal{C}_{k-1}$.  We thus have that a formal map encodes  a generalisation of the notion of a $\pi$-system as introduced by the Turaev in \cite{turaev:hqft1}.

(ii)  The category of crossed complexes is monoidal closed. The crossed complex  of morphisms, $\mathcal{CRS}(\mathcal{C}, \mathcal{D})$, mentioned above, is the value of the functor $\mathcal{CRS}(\mathcal{C}, -)$ on a crossed complex $\mathcal{D}$ and that functor is right adjoint to a tensor product $ -\otimes \mathcal{C}$.  Thus in the description of the nerve of $\mathcal{CRS}(\mathcal{C}, \mathcal{D})$, we have that it  is given in dimension $n$ by
\begin{eqnarray*}Ner(\mathcal{CRS}(\mathcal{C}, \mathcal{D}))_n &=& \mathbf{Crs}(\pi(n), \mathcal{CRS}(\mathcal{C}, \mathcal{D}))\\
							    & \cong & \mathbf{Crs}(\pi(n) \otimes \mathcal{C}, \mathcal{D}).
\end{eqnarray*}
The tensor product is given in terms of generators and relations and so explicit descriptions of maps from a tensor product to another crossed complex are fairly easy to specify.  In particular  for $X = |K|$, a polyhedral space (typically a triangulated manifold), the characteristic maps correspond to 0-simplices in $\mathcal{S}(K, Ner(\mathcal{C}))$, the values of a given HQFT, $\tau$, on a specific $(X,g)$ depend, up to isomorphism,  only on the homotopy class of $g$ and a homotopy is a 1-simplex in this simplicial set,  $\mathcal{S}(K, Ner(\mathcal{C}))$. It is thus given by a morphism of crossed complexes, $\pi(1) \otimes \pi(K) \to \mathcal{C}$, i.e., a homotopy of crossed complex morphisms.  One can specify the homotopy relation, and, if need be, even the homotopy itself, combinatorially by stating where generating cells get sent.

This sort of analysis can also be given at the purely simplicial level leading to a homotopy of simplicial maps from $K$ to $Ner( \mathcal{C})$.  The advantage of the crossed complex approach is that we can replace $\pi K$, defined simplicially, by $\pi |K|$ defined via any regular CW-decomposition of
 $|K|$, which will be completely independent of the choice of order on the vertices of the underlying simplicial complex and may be much smaller and nearer to the `geometry'.  This follows from the general `yoga' of crossed complexes and their relation with simplicial sets (e.g., the crossed complex version of the Eilenberg-Zilber theorem).  The simplicial approach, however, also has its advantages, in particular because of the similarity with lattice based models in TQFTs and the explicit combinatorial / geometric gadgetry available. 


\section{Formal HQFTs}
The notion of a simplicial formal $\mathcal{C}$-map and the corresponding formal $\mathcal{C}$-cobordisms allow us to extend the definition of formal HQFT that we introduced in \cite{TPVT:hqft1} to all dimensions and a general crossed complex, $\mathcal{C}$.
\subsection{Formal structures of formal $\mathcal{C}$-maps}

Before we can give the definition of a formal HQFT, we need to describe some of the constructions we will use.

Supposing that we are working with $d$-dimensional manifolds,  we will need to consider these together with the corresponding cobordisms.  First we note that if $K$ is the empty simplicial complex, for instance, triangulating the empty $d$-dimensional manifold, then there is a unique formal $\mathcal{C}$-map defined on $K$. Next if $\lambda_i : K_i \to Ner(\mathcal{C})$ for $i = 1,2$ are two formal $\mathcal{C}$-maps, then they naturally give a formal $\mathcal{C}$-map, $\lambda_1\sqcup \lambda_2 : K_1\sqcup K_2 \to Ner(\mathcal{C})$, given by the universal property of the coproduct, and unique up to isomorphism given a choice of that coproduct in the usual way. We say this is the \emph{sum} of the  two formal maps.  A certain amount of care needs to be taken in the usual way as different ordered representations give different decompositions.

We will say that a formal $\mathcal{C}$-map, $\lambda : K\to Ner(\mathcal{C})$, is \emph{connected} if the underlying domain, $K$, is a connected simplicial complex.   Given a general formal $\mathcal{C}$-map $\lambda: K \to Ner(\mathcal{C})$, and an ordered decomposition  of $K$ as a disjoint union of its connected components, then, naturally, we get a decomposition of $\lambda$ as a sum of connected formal maps.

If $\mathbf{\Lambda} : \lambda_0\to \lambda_1$ and $\mathbf{\Gamma} :\lambda_1\to \lambda_2$ are two formal  $\mathcal{C}$-cobordisms (with suitable triangulating simplicial complexes subsummed in the notation), then we can construct a composite formal  $\mathcal{C}$-cobordisms  in the obvious way, which we will denote by $\mathbf{\Lambda}\#_{\lambda_1}\mathbf{\Gamma}$.  (If extra structure (e.g., differential manifold structures) is being considered on the manifolds, it will be necessary to use cobordisms with a collar neighbourhood of the boundaries to ensure composition works at the deeper level. Ways of handling this are well known for TQFTs and cause no real problem.)
\subsection{The definition} \label{deffhqft} 
Fix, as before, a crossed complex, $\mathcal{C}$, and also fix a ground field, $\mathbb{K}$.

A \emph{(simplicial) formal HQFT} with background $\mathcal{C}$ assigns\begin{itemize}
\item to each connected (simplicial) formal $\mathcal{C}$-map, $\lambda$,   a $\mathbb{K}$-vector space $\tau(\lambda)$, and by extension, to each formal $\mathcal{C}$-map on a $d$-manifold $X$, given by a list $\lambda = \{\lambda_i\, | \, i \in I\}$ of formal connected $\mathcal{C}$-maps, a tensor product $$\tau(\lambda) = \bigotimes_{i\in I}\tau(\lambda_i);$$
\item to any equivalence class of (simplicial) formal $\mathcal{C}$-cobordisms, $(M,\mathbf{\Lambda})$ between $(X_0,\lambda_0)$ and $(X_1,\lambda_1)$, a $K$-linear transformation 
$$\tau(\mathbf{\Lambda}) : \tau(\lambda_0) \to \tau(\lambda_1),$$
\end{itemize} 
These assignments are to satisfy the following axioms:
\begin{enumerate}[(i)]
\item Disjoint union of formal $\mathcal{C}$-maps corresponds to tensor product of the corresponding vector spaces via specified isomorphisms:
$$\tau(\lambda_0\sqcup \lambda_1) \stackrel{\cong}{\to}\tau(\lambda_0)\otimes \tau(\lambda_1),$$
$$\tau(\emptyset)\stackrel{\cong}{\to}\mathbb{K}$$
for the ground field $\mathbb{K}$, so that a) the diagram of specified isomorphisms
$$\xymatrix{\tau(\lambda)\ar[r]^\cong\ar[d]_\cong &\tau(\lambda\sqcup \emptyset)\ar[d]^\cong \\
\tau(\lambda)\otimes \mathbb{K}&\tau(\lambda)\otimes \tau(\emptyset)\ar[l]^\cong}$$
for $\lambda \to \emptyset\sqcup \lambda$, commutes and similarly for $\lambda\to \lambda\sqcup\emptyset$, and b) the assignments are compatible with the associativity isomorphisms for $\sqcup$ and $\otimes$, so that $\tau$ satisfies the usual axioms for a symmetric monoidal functor.
 \item For formal $\mathcal{C}$-cobordisms
 $$\mathbf{\Lambda} : \lambda_0\to \lambda_1, \quad \mathbf{\Gamma} : \lambda_1\to \lambda_2$$ with composite $\mathbf{\Lambda}\#_{\lambda_1}\mathbf{\Gamma}$, we have 
 $$\tau(\mathbf{\Lambda}\#_{\lambda_1}\mathbf{\Gamma}) = \tau(\mathbf{\Gamma})\tau(\mathbf{\Lambda}) : \tau(\lambda_0) \to  \tau(\lambda_2).$$
 \item For the identity formal $\mathcal{C}$-cobordism on  $\lambda$, 
 $$\tau(1_\lambda) = 1_{\tau(\lambda)}.$$
 \item Interaction of cobordisms and disjoint union is transformed correctly by $\tau$, i.e., given formal $\mathcal{C}$-cobordisms
 $$\mathbf{\Lambda} : \lambda_0\to \lambda_1, \quad \mathbf{\Gamma} : \gamma_0\to \gamma_1,$$ 
 the following diagram
$$ \xymatrix{\tau(\lambda_0\sqcup\gamma_0)\ar[r]^\cong\ar[d]_ {\tau(\mathbf{\Lambda}\sqcup\mathbf{\Gamma})}&\tau(\lambda_0)\otimes\tau(\gamma_0)\ar[d]^{\tau(\mathbf{\Lambda})\otimes\tau(\mathbf{\Gamma})}\\
 \tau(\lambda_1\sqcup \gamma_1)\ar[r]_\cong&\tau(\lambda_1)\otimes \tau(\gamma_1)}$$
 commutes, compatibly with the associativity structure.
\end{enumerate}
\textsc{Remark.}
Replacing the `simplicial' by `cellular' etc. gives a wider definition of formal HQFT and, of course, this has an advantage of allowing smaller calculations for manifolds as there are fewer cells in a CW-decomposition than simplices in a triangulation, in general.
\subsection{The category of formal $\mathcal{C}$-maps}
One idea of a homotopy quantum field theory is that it is a representation of the monoidal category of $B$-cobordisms. 
This was made explicit by Rodrigues, \cite{rodrigues}, who proved that the category, $\mathbf{HCobord}(d,B)$, of $d$-dimensional $B$-manifolds and (homotopy) $B$-cobordisms is a symmetric monoidal category. (A similar observation had been made by Brightwell and Turner \cite{B&T} on the low dimensional case of the homotopy surface category, linked to constructions of Segal, Tillmann and others.)  With that interpretation,  a $(d+1)$-HQFT is a means of studying $\mathbf{HCobord}(d,B)$ via a representation, i.e., a monoidal functor from $\mathbf{HCobord}(d,B)$ to the category of vector space over some field or, more generally, to any well understood and nicely behaved symmetric monoidal category.

Given the motivation of these papers, it seems clear that there should be a symmetric monoidal category of   (simplicial) formal $\mathcal{C}$-maps so that a formal HQFT with $\mathcal{C}$ as base was a symmetric monoidal functor from it to $Vect$.  This is more or less clear but needs a little care in the setting up.

 We let $\mathbf{FHCobord}(d,\mathcal{C})$ have the following claimed categorical structure:
\begin{itemize}
\item its objects are oriented $d$-dimensional manifolds $X$, each together with a triangulation $\mathbf{T}$ and a formal $\mathcal{C}$-map $\lambda : T \to Ner(\mathcal{C})$;
\item its morphisms are equivalence classes of formal $\mathcal{C}$-cobordisms between such formal $\mathcal{C}$-maps;
\item its composition is given by gluing of cobordisms in the obvious way;
\item for a given $(X,\mathbf{T},\lambda)$, the corresponding identity is the equivalence class of the cylinder cobordism on $X\times I$ with triangulation and $\mathcal{C}$-coloring as considered earlier;  
\item the monoidal category structure is given by `coproduct over $Ner(\mathcal{C})$', that is, given $(X_i,\mathbf{T}_i,\lambda_i)$ for $i = 1,2$, we take the disjoint union of the manifolds $X_1\sqcup X_2$ with the obvious induced triangulation giving a simplicial complex $T_1\sqcup T_2$ and then use the universal property for coproduct / disjoint union to give the map to $Ner(\mathcal{C})$;
\item the unit of the monoidal structure is the empty formal $\mathcal{C}$-map.
\end{itemize}
\begin{theorem}
The above definition makes $\mathbf{FHCobord}(d,\mathcal{C})$ into a symmetric monoidal category.
\end{theorem}
\textsc{Proof.}
Most of this is routine as similar arguments are well represented in the literature on TQFTs. One point of note is that the category structure, and in particular, the identities of that structure, is where it becomes necessary to work with equivalence classes of cobordisms, and not just with the formal $\mathcal{C}$-cobordisms themselves. The sort of argument is well known. Attaching a cylinder to an incoming or outgoing boundary of a cobordism changes the cobordism, but does keep within the equivalence class. \hfill $\square$

\,

The following is now an obvious reformulation / corollary of this result.  In the case that $\mathcal{C}$ is a crossed complex with an abelian group $A$ in dimension 2 and trivial groups everywhere else, the formal HQFTs on $\mathcal{C}$ are exactly the HQFTs with background $K(A,2)$ considered by Brightwell and Turner in \cite{B&T} and so this result extends the corresponding observation in their work.
\begin{theorem}
A (simplicial) formal HQFT, $\tau$, with background $\mathcal{C}$ corresponds to a representation
$$\tau : \mathbf{FHCobord}(d,\mathcal{C})\to Vect.$$\hfill $\square$
\end{theorem}
No essential role is played by any simplicial hypothesis here and so one should expect  similar result for theories based on cellular or handle decompositions on the one hand and ones in which $\mathcal{C}$ is replaced by  simplicial group on the other.

It is worth noting that equivalent $d$-dimensional formal $\mathcal{C}$-maps on a manifold give isomorphic objects in $\mathbf{FHCobord}(d,\mathcal{C})$, so effectively are independent of the decomposition used.
\subsection{Operations on Formal HQFTs}
As a formal HQFT is a monoidal functor to $Vect$, (but this could be replaced by any other suitably nice additive monoidal category), there are some obvious operations that can be performed on them, just as in the non-formal case.

If $\tau$ and $\tau^\prime$ are two FHQFTs with the same background, $\mathcal{C}$, then set $(\tau \oplus \tau^\prime)(\lambda) = \tau(\lambda) \oplus \tau^\prime(\lambda)$ and similarly for the $\mathcal{C}$-cobordisms, to get their direct sum.  The tensor product $(\tau \otimes \tau^\prime)$ is defined similarly by using the tensor in the image category.  There is also a dual given by noting that any formal $\mathcal{C}$-cobordism $(M,\mathbf{\Lambda})$ can be reversed to get a cobordism in the opposite direction, which can be considered with the reverse orientation (see \cite{turaev:hqft1} and \cite{rodrigues} for this in the non-formal setting).

We define a category $FQ_{d+1}(\mathcal{C})$ to be that of all formal $(d+1)$-HQFTs with given background $\mathcal{C}$ and the natural monoidal transformations between them.

Now suppose $p : \mathcal{D}\to \mathcal{C}$ is a morphism of crossed complexes. This induces a strict monoidal functor
$$\mathbf{FHCobord}(d,\mathcal{D})\to \mathbf{FHCobord}(d,\mathcal{C})$$
by composition, $\lambda$ being sent to $Ner(p)\lambda$.  As $FQ_{d+1}(\mathcal{C})$ is just the category of monoidal `representations' of $\mathbf{FHCobord}(d,\mathcal{C})$, this clearly induces a functor
$$p^* : FQ_{d+1}(\mathcal{C})\to FQ_{d+1}(\mathcal{D}).$$
One would expect that $p^*$ might have left and right adjoints analogously to the usual setting of representations.  The usual method would be to use Kan extension type formulae and to see if the necessary limits or colimits exist.  Of course, $Vect $ is the category of \emph{finite dimensional} vector spaces, so we should expect not to be able to find adjoints for all $p$, even with finite crossed complexes as both domain and codomain, or perhaps more exactly, the proof that  an adjoint exists (or not) is likely to depend on properties of $p$.  In fact we will assume shortly that $p$ is a fibration of crossed complexes or more exactly that $Ner(p)$ is a Kan fibration and that $Ker(p)$ is a finite crossed complex.  These conditions are almost certainly stronger than necessary, but are useful here as the proofs of existence are fairly easy  \emph{and} they extend known results, for instance, results in \cite{turaev:hqft1} on push-forward and transfer.  We also restrict attention to the left adjoint to $p^*$ corresponding to a right Kan extension.  The relevant formula in more or less standard notation is: for a formal $\mathcal{C}$-map, $\lambda$,
$$R_p(\tau)(\lambda) = Colim((p^*\downarrow \lambda) \stackrel{\delta}{\to} \mathbf{FHCobord}(d,\mathcal{D}) \stackrel{\tau}{\to} Vect).$$
An object of comma category $(p^*\downarrow \lambda) $ consists of a pair  $(\mu, \Gamma)$ where $\mu :T \to Ner(\mathcal{D})$ is a formal $\mathcal{D}$-map and $\Gamma : p(\mu) \to \lambda$ is a formal $\mathcal{C}$-cobordism. The functor $\delta$ sends $(\mu, \Gamma)$ to $\mu$ and $\tau$ is the formal $\mathcal{D}$-HQFT being considered.

The morphisms of $(p^*\downarrow \lambda) $ from $(\mu, \Gamma)$ to $(\mu^\prime, \Gamma^\prime)$, say, are the formal $\mathcal{D}$-cobordisms, $\mathbf{\Phi}:\mu \to \mu^\prime$, such that $\Gamma^\prime \# p(\mathbf{\Phi}) = \Gamma$.   Such a morphism gives a linear transformation
$$\tau(\mathbf{\Phi}) : \tau(\mu)\to \tau(\mu^\prime).$$ 
Forgetting the finite dimensionality for the moment, the corresponding colimit could be constructed in the usual way as a quotient of a coproduct over all $(\mu, \Gamma)$ of the various $\tau(\mu)$. We can represent elements of the colimit
by symbols $(\mu, \Gamma)\otimes x$ with $x\in \tau(\mu)$, where the equivalence relation determining the quotient is generated by all
$$(\mu,\Gamma^\prime\#p(\mathbf{\Phi}))\otimes x \equiv (\mu^\prime,\Gamma^\prime)\otimes \tau(\mathbf{\Phi})(x).$$ 
If $Ner(p)$ is a fibration, as assumed, then any $\lambda : T \to Ner(\mathcal{C})$ lifts (non-uniquely) to some $\mu : T \to Ner(\mathcal{D})$, similarly for any $\Gamma$.  The number of lifts is finite, as $Ker(p)$ is assumed to be a finite crossed complex. (We will say that $p$ has \emph{finite fibre} in this case.) 
\begin{proposition}
If $p$ is a fibration with finite fibre, then $R_p(\tau)$ exists for any $\tau$.
\end{proposition}
\textsc{Proof.} Any $(\mu, \Gamma)\otimes x$ is equivalent to one in which $p(\mu) = \lambda$ and $\Gamma$ is the identity $\mathcal{C}$-cobordism, hence the colimit is a quotient of a finite direct sum of finite dimensional spaces, so exists in $Vect$.\hfill $\square$
\begin{theorem}
If $p$ is a fibration with finite fibre, then $p^*$ has a left adjoint
$$p_* : FQ_{d+1}(\mathcal{D})\to  FQ_{d+1}(\mathcal{C}).$$
\end{theorem}
\textsc{Proof.} We, of course, define $p_*(\tau) = R_p(\tau)$, defined as above.  It is fairly routine to check that it is a formal $\mathcal{C}$-HQFT, but it may help if we note that if $\lambda \sqcup \lambda^\prime$ is a sum of formal $\mathcal{C}$-maps on a disjoint union, then 
$$(p^*\downarrow \lambda \sqcup \lambda^\prime) \simeq (p^*\downarrow \lambda)\sqcup (p^*\downarrow \lambda^\prime),$$
and that, if $\mathbf{\Lambda} : \lambda \to \lambda^\prime$ is a formal $\mathcal{C}$-cobordism, then
$$p_*(\tau)(\mathbf{\Lambda})((\mu,\Gamma)\otimes x) = ((\mu, \mathbf{\Lambda}\#\Gamma)\otimes x).$$\hfill$ \square$

\textsc{Example.}
Suppose $d = 1$ and $p : \mathcal{D}\to \mathcal{C}$ is a fibration with finite fibre.  The results of the first part of this series of papers, \cite{TPVT:hqft1}, gave a classification of formal $\mathcal{C}$-HQFTs in terms of crossed $\mathcal{C}$-algebras. (We will not recall these here as they will only be used in this example.)  There, results on pulling back and `pushing forward' such crossed algebras were given, extending results of Turaev, \cite{turaev:hqft1}.  It is easily checked thar, if $L_\tau$ is the crossed $\mathcal{C}$-algebra corresponding to $\tau \in FQ_2(\mathcal{C})$, then, in the notation of that earlier paper, $p^*(L_\tau)$ corresponds to $p^*(\tau)\in FQ_2(\mathcal{D})$ in the sense of our above construction.  It can also easily be checked that the corresponding statement for the two meanings of $p_*$ also holds.
\section{Formal $\mathcal{C}$-maps and combinatorial $\mathcal{C}$-bundles}
This section might have as a subtitle: ``What does the classifying space of a crossed complex classify?''  

If we return to the initial case of a finite group $G$ (so $\mathcal{C}$ would be $G$ in dimension 1 and trivial otherwise), then $|Ner(\mathcal{C})|$ was $BG$, the classifying space of $G$, and it is `classical' that a map $g: X \to BG$ corresponds to an induced principal $G$-bundle on $X$, so giving geometric significance to the characteristic maps of $B$-manifolds and cobordisms when $B = BG$.  The suggested subtitle asks if there is a similar interpretation for a general $\mathcal{C}$ and for formal $\mathcal{C}$-maps.  In the remainder of this paper we will examine some of the approaches to this, summarising the theory from both a crossed complex and a simplicial viewpoint and, hopefully, identifying where further clarification is needed.

In \cite{Attal:2002}, Attal  gives a combinatorial version of non-abelian gerbes with connection and curvature.  This has a definite similarity to the formulation of formal $\mathcal{C}$-maps as given in  \cite{TPVT:hqft1}, (so here $\mathcal{C}$ is a crossed module).  Is there a similar combinatorial `gerbe-style' interpretation of our more general formal $\mathcal{C}$-maps?  The answer is ``Yes, but ...''.  The hesitation is due to there being several versions of partial answers.  We will look at two, one using higher dimensional category theory, the other being a simplicial combinatorial approach.  Each handles some parts of the answer well and yet fails to deliver the whole picture in full generality.  The partial answers however are already significant and deserve an `airing', especially as in some sense they use classical ideas that have perhaps slipped from being `centre stage'.

First an aside and a reformulation:  we have used triangulations of manifolds throughout this paper so far, but another approach is possible and has advantages for the geometric (and physical) interpretations.  Any triangulation of a space $X$ yields an open cover of $X$ by the open stars of vertices of the triangulating simplicial complex. (This can be found in many classical algebraic topology textbooks, usually with regard to \v{C}ech homology and cohomology.) Conversely any open cover of a space has a nerve (in the sense of \v{C}ech) and if one takes the open star cover, as above, its nerve is the same as the simplicial complex used in the triangulation.  Any open cover can be refined to an open star cover of some subdivided triangulation, so open covers should be just as good as triangulations in formulating notions of formal maps on manifolds (or more general spaces).  We therefore give a version of formal $\mathcal{C}$-map relative, not to a triangulation, but to an open cover.  This will make the link with stacks, gerbes and 2-bundles almost immediate.

Just as we needed to order the vertices of the simplicial complexes used earlier, we need here to order the indexing set for the open sets in a given cover $\mathcal{U}$.  An open cover, $\mathcal{U}$, on $X$ together with a total order on its indexing set will be called an \emph{ordered open cover of} $X$.  We can use the order to turn the simplicial complex $N(\mathcal{U})$, which is the nerve of the cover, into a simplicial set.

\textsc{Definition.} Given an ordered open cover $\mathcal{U}$ of $X$ and a crossed complex $\mathcal{C}$, a \textit{(simplicial) formal $\mathcal{C}$-map on $X$ relative to $\mathcal{U}$} is a simplicial map 
$$\lambda : N(\mathcal{U}) \to Ner(\mathcal{C}).$$

\medskip

We start dismantling this definition: suppose $U_i\cap U_j \neq \emptyset $ with $i < j$, then $\langle U_i,U_j\rangle$ is a 1-simplex of $N(\mathcal{U})$, and $\lambda$ assigns some element $\lambda_{ij} \in C_1$ to this.  Likewise if $U_i\cap U_j \cap U_k\neq \emptyset $ with $i < j< k$, then one gets a 2-simplex $\sigma = \langle U_i,U_j,U_k\rangle$ and a corresponding $\lambda(\sigma) \in C_2$ such that $\partial \lambda(\sigma) = \lambda_{02}\lambda_{01}^{-1}\lambda_{12}^{-1}.$ As before these elements and the corresponding conditions continue to higher dimensions.

A complete interpretation of this would require a digression to develop more of the theory of crossed complexes than we have available here so we will restrict to the crossed module case. We know that crossed modules correspond to internal group objects in the category of groupoids and thus to strict 2-groups.  Given a crossed module, $\mathcal{C}$, we can form its associated 2-group $G(\mathcal{C})$.  Referring to Baez and Schreiber, \cite{baez-schreiber} or, more briefly, to Baez's notes for his talks at this conference, we can formulate a notion of principal $G(\mathcal{C})$-2-bundle, with local trivialisation over the open cover, $\mathcal{U}$.  Noting that as $\mathcal{C}$ is a discrete crossed module, the transitions will be constant on the connected components of each intersection, we get using \cite{baez-schreiber}, Prop. 2.3:
\begin{proposition}
If $\mathcal{C}$ is a crossed module with associated strict 2-group $G(\mathcal{C})$, then a formal $\mathcal{C}$-map on $X$ relative to an (nice) open cover $\mathcal{U}$ corresponds to a principal $G(\mathcal{C})$-2-bundle with local trivialisation over the open cover amd hence to a non-abelian gerbe on $X$.
\end{proposition}
\textsc{Proof.} The `niceness' condition is a technical condition\footnote{For manifolds any open cover can be refined to a nice one} that ensures that all the intersections are contractible and hence connected, so the transitions, etc., of a $G(\mathcal{C})$-2-bundle will be constant on each such intersection, hence can be thought of as assigning a value to the corresponding simplex of the nerve.  \hfill $\square$

\medskip

\textsc{Remark.} The treatment of 2-stacks  and 2-gerbes by Breen, \cite{breen1994}, already contains the same formulae, whilst Duskin's \cite{Duskin1989} and Street's \cite{Street} provide  treatments of higher dimensional descent, again containing a discussion of these ideas and linking them with a simplicial treatment.  An application of these to TQFTs in a very closely related context to that of our formal HQFTs can be found in \cite{tqft1,tqft2}.

This proposition gives a partial answer to our problem; `partial' because (i) it requires $\mathcal{C}$ to be a crossed module not a general reduced crossed complex, and (ii) it does not handle formal $\mathcal{C}$-cobordisms or equivalence.  The first of these points will be partially addressed in a forthcoming series of papers by Brown, Glazebrook and the author, \cite{RBJGTP}, which will also look at aspects of \emph{smooth} crossed complexes as well as the discrete case and may thus allow extensions of the theory here to other Lie crossed complexes as formal backgrounds.  It will also look at $\mathcal{C}$-bundles in more generality and their links with cohomology with coefficients in $\mathcal{C}$.  The second point can be partially answered by examining the corresponding discussion in  \cite{tqft1,tqft2} and adapting the treatment according to the different context.  This is clearly do-able, but has been put off to a later paper.  A full answer will require a much fuller discussion of the  theory of crossed complexes and their relatives than there was space, or time, for here.  That it can, and probably should, be done is indicated by the second set of partial answers using a combinatorial and thus simplicial approach.

The full simplicial geometric interpretation of formal $\mathcal{C}$-maps is still some way off, however the main lines are clear.  Given a formal $\mathcal{C}$-map
$$\lambda : T \to Ner(\mathcal{C})$$
on a triangulation of a manifold $X$ or a cobordism $M$, we saw earlier that $\mathcal{C}$ corresponded to the Moore complex of a simplicial group, which we will also denote $\mathcal{C}$ here, and $Ner(\mathcal{C})$ could be realised as $\overline{W}(\mathcal{C})$. 

Classically (cf. Curtis, \cite{curtis} \S 6) to any simplicial group $G$, this classifying space $\overline{W}(G)$ comes together with a total complex $W(G)$ and a projection
$$p : W(G) \to \overline{W}(G)$$
yielding a classifying $G$-bundle.  That theory of simplicial fibre bundles also gives us that any simplicial fibre bundle has a description as a \emph{regular twisted Cartesian product}, relative to a twisting function $t : B \to G$ or in general into $aut(F)$ where $F$ is the fibre.  Here $B$ is the base. The twisting function $t$ corresponds exactly to a simplicial map $f_t : B \to \overline{W}(G)$. The formulae can be found in the survey by Curtis mentioned earlier.  

The theory of simplicial fibre bundles thus follows a parallel track to the better known topological theory, but with the twisted Cartesian product giving a neat combinatorial way of handling them.  The twisting functions are `the same' as simplicial maps to the classifying space and that suits us well because in our context those are the formal $\mathcal{C}$-maps. We thus suggest a definition:

\textsc{Definition.} If $\mathcal{C}$ is a crossed complex (and also its associated simplicial group), a \emph{combinatorial $\mathcal{C}$-bundle on $X$ relative to the triangulation} $T$ (or \emph{to an open cover} $\mathcal{U}$) will be a principal simplicial $\mathcal{C}$-bundle on $T$ (resp. on $N(\mathcal{U})$).
\begin{proposition}
Any formal $\mathcal{C}$-map on $X$ relative to $T$ (or $\mathcal{U}$) corresponds to a combinatorial $\mathcal{C}$-bundle on $X$ rel. $T$ (or rel. $\mathcal{U}$).
\end{proposition}
\textsc{Proof.} A formal $\mathcal{C}$-map $\lambda$ gives the twisted Cartesian product $T\times_\lambda \mathcal{C}$. \hfill $\square$

The treatment of formal $\mathcal{C}$-cobordisms is clear, and the treatment in \cite{tqft1,tqft2} suggests the way to give a detailed treatment of equivalence.

\textsc{Remark.} Simplicial fibre bundle theory was developed in the late 1960s and early 1970s.  Its development was influenced by its application to geometric problems relating to triangulations and smoothings. There the simplicial groups were not finite, but similarities with current problems in QFTs and the natural way the simplicial theory emerges in HQFTs as a combinatorial  approach suggest that the old results and methods may deserve being disinterred and dusted off with a view to their adaptation to modern problems.

The  combinatorial / simplicial theory gives a moderately complete answer to our problem, but it needs sheafifying and the links between it and the higher dimensional categorical approach need some clarification.  It strengthens the perception that TQFTs and HQFTs are in some sense non-abelian analogues of $K$-theory.  Finally it should be mentioned that all this stack / gerbe / formal map machine is very closely related to Grothendieck's Pursuit of Stacks, which returns us to the problems dear to Ross Street and to one of the recurrent themes of talks at this meeting.


\end{document}